\magnification 1200
\input amstex.tex
\input amsppt.sty

\def\authorhead{ NAZAROV AND TARASOV }
\def\titlehead{  TENSOR PRODUCTS OF YANGIAN MODULES }

\hfuzz 10pt
\NoBlackBoxes


\hsize 5.208in	
\vsize 8.025in	
\hoffset 0pt
\voffset 13truept 


\edef\enddef{\catcode`\noexpand\@=\the\catcode`\@\space}

\catcode`\@=11
\mathsurround 1.6pt

\font\Bbf=cmbx12 scaled 1200


\def\leftheadline{\rlap{\foliofont@\folio}\hfil\ignorespaces\authorhead
 \unskip\hfil}
\def\rightheadline{\hfil\ignorespaces\titlehead\unskip\hfil
 \llap{\foliofont@\folio}}

\headline={\def\\{\unskip\space\ignorespaces}\headlinefont@
 \ifodd\pageno \rightheadline \else \leftheadline\fi}
\footline={\hfil}

\def\makeheadline{\vbox to\z@{\vskip-22.5\p@
 \line{\vbox to8.5\p@{}\the\headline}\vss}\nointerlineskip}
\def\makefootline{\baselineskip24\p@\lineskiplimit\z@\line{\the\footline}}


\def\varitem#1#2{\par\setbox\z@\hbox{#1\enspace}\hangindent\wd\z@
 \hglue-2\parindent\kern\wd\z@\indent\llap{#2\enspace}\ignore}

\def\bitem{\varitem{\space(b)}}

\def\refskips{\frenchspacing \let\csk@p\, \def\,{\allowbreak\csk@p}}
\def\refitem{\varitem{[VW]}} \def\widest#1{\def\refitem{\varitem{#1}}}

\let\bls\baselineskip  \let\ignore\ignorespaces
\def\vsk#1>{\vskip#1\bls} \let\adv\advance 
\def\vv#1>{\vadjust{\vsk#1>}\ignore}
\def\vvn#1>{\vadjust{\nobreak\vsk#1>\nobreak}\ignore}
\def\vvgood{\vadjust{\penalty-500}}

\def\emph#1{{\it #1\/}}

 \let\nt\noindent 
\def\nn#1>{\noalign{\vskip #1pt}} \def\NN#1>{\openup#1pt}
 
\let\Sum\sum \def\sum{\Sum\limits} 
\let\Prod\prod \def\prod{\Prod\limits} \let\Int\int \def\int{\Int\limits}

\let\=\m@th \def\&{.\kern.1em} \def\>{\!\;} \def\:{\!\!\;}

\ifx\plainfootnote\undefined \let\plainfootnote\footnote \fi
\expandafter\ifx\csname amsppt.sty\endcsname\relax
 
\else \fi

\newbox\s@ctb@x
\def\s@ct#1 #2\par{\removelastskip\vsk>
 \vtop{\bf\setbox\s@ctb@x\hbox{#1} \parindent\wd\s@ctb@x
 \ifdim\parindent>0pt\adv\parindent.5em\fi\item{#1}#2\strut}%
 \nointerlineskip\nobreak\vtop{\strut}\nobreak\vsk-.4>\nobreak}

\newbox\t@stb@x
\def\gadv{\global\advance} \def\gad#1{\gadv#1 1} 
\def\l@b@l#1#2{\def\n@@{\csname #2no\endcsname}%
 \if *#1\gad\n@@ \expandafter\xdef\csname
 @#1@#2@\endcsname{\the\Sno.\the\n@@}%
 \else\expandafter\ifx\csname @#1@#2@\endcsname\relax\gad\n@@
 \expandafter\xdef\csname @#1@#2@\endcsname{\the\Sno.\the\n@@}\fi\fi}
\def\l@bel#1#2{\l@b@l{#1}{#2}\?#1@#2?}
\def\?#1?{\csname @#1@\endcsname}
\def\[#1]{\def\n@xt@{\ifx\t@st *\def\n@xt####1{{\setbox\t@stb@x\hbox{\?#1@F?}%
 \ifnum\wd\t@stb@x=0 {\bf???}\else\?#1@F?\fi}}\else
 \def\n@xt{{\setbox\t@stb@x\hbox{\?#1@L?}\ifnum\wd\t@stb@x=0 {\bf???}\else
 \?#1@L?\fi}}\fi\n@xt}\futurelet\t@st\n@xt@}
\def\(#1){{\rm\setbox\t@stb@x\hbox{\?#1@F?}\ifnum\wd\t@stb@x=0 ({\bf???})\else
 (\?#1@F?)\fi}}
\def\dff{\expandafter\d@f} \def\d@f{\expandafter\def}
\def\edff{\expandafter\ed@f} \def\ed@f{\expandafter\edef}

\newcount\Sno \newcount\Lno \newcount\Fno
\def\Section#1{\gad\Sno\Fno=0\Lno=0\s@ct{\the\Sno.} {#1}\par}
 \let\Sect\Section
\def\section#1{\gad\Sno\Fno=0\Lno=0\s@ct{} {#1}\par} \let\sect\section
\def\l@F#1{\l@bel{#1}F} \def\<#1>{\l@b@l{#1}F} \def\l@L#1{\l@bel{#1}L}
\def\Tag#1{\tag\l@F{#1}} \def\Tagg#1{\tag"\llap{\rm(\l@F{#1})}"}
\def\Th#1{Theorem \l@L{#1}} \def\Lm#1{Lemma \l@L{#1}}
\def\Prop#1{Proposition \l@L{#1}}
\def\Cr#1{Corollary \l@L{#1}} \def\Cj#1{Conjecture \l@L{#1}}
 
\def\Proof#1.{\demo{\it Proof #1}}

 \def\setparindent{\edef\Parindent{\the\parindent}}
\def\Appendix{\Sno=64\let\p@r@\z@ 
\def\Section##1{\gad\Sno\Fno=0\Lno=0 \s@ct{}
 \hskip\p@r@ Appendix \char\the\Sno
 \if *##1\relax\else {.\enspace##1}\fi\par} \let\Sect\Section
\def\section##1{\gad\Sno\Fno=0\Lno=0 \s@ct{} \hskip\p@r@ Appendix%
 \if *##1\relax\else {.\enspace##1}\fi\par} \let\sect\section
\def\l@b@l##1##2{\def\n@@{\csname ##2no\endcsname}%
 \if *##1\gad\n@@
\expandafter\xdef\csname @##1@##2@\endcsname{\char\the\Sno.\the\n@@}%
\else\expandafter\ifx\csname @##1@##2@\endcsname\relax\gad\n@@
 \expandafter\xdef\csname @##1@##2@\endcsname{\char\the\Sno.\the\n@@}\fi\fi}}

\let\logo@\relax

\def\>{\RIfM@\mskip.6\thinmuskip\relax\else\kern.1em\fi}
\def\]{\RIfM@\mskip-.6\thinmuskip\relax\else\kern-.1em\fi}
\let\ts\> \let\ns\]


\def\A{\text{\rm A}}
\def\al{\alpha}
\let\alb\allowbreak
\def\astr{{\raise1.5pt\hbox{\kern-6.5pt$^{\,\ast}$}}}

\def\D{\Cal D}
\def\be{\beta}
\def\bx{{\boxed{\phantom{\square}}\kern-.4pt}}

\def\CC{{\Bbb C}}
\def\chip{\chi^{\ts\prime}}
\def\cp{c^{\ts\prime}}

\def\crc{{\raise.24ex\hbox{$\sssize\kern.1em\circ\kern.1em$}}}
\def\CS{\CC\hskip-.5pt\cdot\hskip-1pt S}

\def\db{d^{\,\prime}}
\def\De{\Delta}
\def\de{\delta}
\def\Det{\widetilde\Delta}

\def\End{\mathop{\roman{End}\>}\nolimits}
\def\enddemos{{$\quad\square$\enddemo}}
\def\ep{\varepsilon}

\def\Fv{F^{\hskip.5pt\vee}}

\def\ga{\gamma}
\def\ge{\geqslant}
\def\glN{\frak{gl}_N}
\def\Gv{G^{\hskip.5pt\vee}}

\def\id{\roman{id}\ts}

\def\io{\iota}

\def\L{\Cal L}
\def\la{\lambda}
\def\lap{\lambda^\ast}
\def\lcd{{\ts,\hskip.95pt\ldots\ts,\ts}}
\def\le{\leqslant}
\def\lox{\lsym\ox}
\def\lsym#1{#1\alb\ldots\relax#1\alb}

\def\mi{{\raise.5pt\hbox{-}}}
\def\mup{\mu^\ast}
\def\mv{\kern127pt}
\def\mw{\kern-81pt}

\def\Om{\Omega}
\def\om{\omega}
\def\omf{{\om^\flat}}
\def\omn{{\om^\natural}}
\def\omp{{\om^\ast}}
\def\ompr{\om^{\ts\prime}}
\def\or{(0\lcd\ns0)}
\def\ox{\otimes}
\def\oxp{\mathbin{\widetilde{\otimes}}}

\def\ph{\varphi}
\def\ps{{\hskip-.5pt\psi}}

\def\Q{\Cal Q}

\def\R{\Cal R}
\def\RR{{\Bbb R}}

\def\si{\sigma}
\def\sgn{\operatorname{sgn}}

\def\tr{\mathop{\roman{tr}\>}\nolimits}
\def\Tti{\widehat T}

\def\U{U}
\def\Ux#1{U^{\ox #1}}

\def\YN{\operatorname{Y}(\glN)}

\def\ZZ{{\Bbb Z}}

\enddef


\line{\Bbf On irreducibility of tensor products of Yangian\hfill}
\medskip
\line{\Bbf modules associated with skew Young diagrams\hfill}
\bigskip\bigskip\bigskip\bigskip
\line{MAXIM NAZAROV$^{1}$ and VITALY TARASOV$^{2}$\hfill}
\medskip
\line{$^{1}$\it
Department of Mathematics, University of York, York YO1 5DD, England\hfill}
\line{\phantom{$^{1}$}\it
E-mail\ts: mln1\@\>york.ac.uk\hfill}
\smallskip
\line{$^{2}$\it
Steklov Mathematical Institute, Fontanka 27, St.\,Petersburg 191011,
Russia\hfill}
\line{\phantom{$^{2}$}\it
E-mail\ts: vt\@\>pdmi.ras.ru\hfill}
\bigskip\bigskip\bigskip\bigskip 


\line{\bf Introduction\hfill}
\medskip\smallskip\nt
In this article we continue our study \text{[\ts NT1\ts,\ts NT2\ts]}
of the finite-dimensional
modules over the Yangian $\YN$ of the general linear Lie algebra $\glN$.
The Yangian $\YN$ is a canonical deformation of the
enveloping algebra $\operatorname{U}(\glN[u])$
in the class of Hopf algebras [D1]. The unital associative
algebra $\YN$ has an infinite
family of generators $T_{ij}^{(s)}$ where $s=1,2,\ts\ldots\ts$ and
$i\ts,\ns j=1\lcd N$. In Section 4 we recall the defintion of the
Hopf algebra $\YN$ in terms of the generating~series
$$
T_{ij}(u)=\de_{ij}+T_{ij}^{(1)}u^{-\ns1}+T_{ij}^{(2)}u^{-\ns2}+\ldots\,\ts.
\tag{0.1}
$$
The classification of the irreducible finite-dimensional $\YN\]$-modules
has been given by V.\,Drinfeld [D2]
by generalizing the results of the second author
\text{[\ts T1\ts,\ts T2\ts]\ts.}
But the structure of these modules still needs better
understanding. For instance, the dimensions
of these modules are not explicitly known in general.

There is a distinguished family of irreducible finite-dimensional
$\YN\]$-modules [\ts C2\ts,O\ts]. We call these modules elementary.
Every elementary $\YN\]$-module is determined by a complex number $z$,
and by two partitions
$\la=(\la_1,\la_2\ts,\,\ldots)$ and $\mu=(\mu_1,\mu_2\ts,\,\ldots)$
such that $\la_i\ge\mu_i$ for any $i=1,2,\ts\ldots\ts$.
This module depends on $\la$ and $\mu$ through
the skew Young diagram  $\om=\la\ts/\mu\ts$, see Section 1.
We denote the corresponding $\YN\]$-module by $V_\om(z)$,
its explicit construction will be given in Section~4.
The vector spaces of the modules $V_\om(z)$ are the same for all
values of $z$, we denote this vector space by $V_\om$.
Suppose that for some positive integer $M$
the numbers of non-zero parts of $\la$ and $\mu$ do not exceed
$M\!+\ns N$ and $M$ respectively. Then the dimension of $V_\om$ equals
the multiplicity of the irreducible $\frak{gl}_M\]$-module of highest
weight $(\mu_1\lcd\mu_M)$ in the restriction to $\frak{gl}_M$
of the the irreducible $\frak{gl}_{M+N}\]$-module of
highest weight $(\la_1\lcd\la_{M+N})$.
In particular, the space $V_\om$ is non-zero if and only if
$\lap_i-\mup_i\le N$ for any index $i=1,2,\ts\ldots\ts$.
Here $\lap=(\lap_1,\lap_2\ts,\,\ldots)$ and $\mup=(\mup_1,\mup_2\ts,\,\ldots)$
are the partitions conjugate to $\la$ and $\mu$ respectively.

For any formal series $f(u)\in \CC\ts[[u^{-1}]]$ with the
leading term 1, the assignment $T_{ij}(u)\mapsto f(u)\cdot T_{ij}(u)$
determines an automorphism of the algebra $\YN$. Let us denote by $\al_f$ this
automorphism. Further, there is a canonical chain of algebras
$\operatorname{Y}(\frak{gl}_1)\subset\ldots\subset
\operatorname{Y}(\frak{gl}_N)$.
The elementary modules are distinguished amongst all irreducible
finite-dimensional $\YN\]$-modules $V$  by the following fact.
Consider the commutative subalgebra in $\YN$ generated by the centres
of all algebras in that chain. This subalgebra is maximal commutative
[\ts C2\ts,\ts NO\ts]\>.
The action of this subalgebra in $V$
is semi-simple if and only if the module $V$
is obtained by pulling back through some
automorphism $\al_f$ from the tensor product
$$
V_{\om_1}(z_1)\ox\ldots\ox V_{\om_k}(z_k)
\tag{0.2}
$$
of the elementary $\YN\]$-modules,
for some skew Young diagrams $\om_1\lcd\om_k$ and some complex numbers
$z_1\lcd z_k$ such that $z_i-z_j\notin\ZZ$ when $i\neq j$.
This fact was
conjectured by I.\,Cherednik, and proved by him in [C2]
under certain extra conditions
on the module $V$. In full generality this fact has been proved in
\text{ [NT1]\ts.}

Furtermore, if $z_i-z_j\notin\ZZ$ for all $i\neq j$,
the tensor product of $\YN\]$-modules (0.2) is always irreducible.
Given the skew Young diagrams $\om_1\lcd\om_k$ it is natural to ask,
for which values of $z_1\lcd z_k\in\CC$ the $\YN\]$-module (0.2) is
reducible. In [NT2] we answered this question
in the case when each of the the diagrams $\om_1\lcd\om_k$ has a rectangular
shape. In that particular case, our answer confirmed the following general
conjecture. For any indices $i$ and $j$ consider the tensor product
$V_{\om_i}(z_i)\ox V_{\om_j}(z_j)$. Consider also the tensor product
$V_{\om_i}(z_i)\oxp V_{\om_j}(z_j)$ obtained via the opposite comultiplication
on $\YN$. If $z_i-z_j\notin\ZZ$,
these two tensor products are irreducible and
equivalent. Hence there is an invertible intertwining operator
$$
R_{ij}(z_i,\ns z_j):\ts
V_{\om_i}(z_i)\ox V_{\om_j}(z_j)\longrightarrow
V_{\om_i}(z_i)\oxp V_{\om_j}(z_j)\,,
$$
unique up to a non-zero multiplier from $\CC$. These multipliers
can be chosen so that $R_{ij}(z_i,\ns z_j)$ depends on $z_i$ and $z_j$
as a rational function of $z_i-z_j$, see Section~4. Now for any
$z_i-z_j\in\CC$ denote by $I_{ij}$ the leading coefficient
of the expansion of the rational function
$R_{ij}(z_i\ns+\ns u\ts,\ns z_j)$ into Laurent series in $u$
near the origin $u=0$.
This coefficient is always an interwining operator between
the two tensor products.
If $z_i-z_j\notin\ZZ$, we simply have $I_{ij}=R_{ij}(z_i,\ns z_j)$.
It has been conjectured that
the $\YN\]$-module (0.2) is irreducible if and only if for all $i<j$
the operators $I_{ij}$ are invertible\ts; see for instance [CP].
We prove this conjecture in the present article, see Theorem 4.8.

This theorem implies, in particular, that the $\YN\]$-module (0.2)
is irreducible if and only if for all $i<j$ the pairwise
tensor products $V_{\om_i}(z_i)\ox V_{\om_j}(z_j)$ are irreducible.
Given the skew Young diagrams $\om_i$ and $\om_j$,
it would be interesting to describe explicitly the set of all differences
$z_i-z_j\in\ZZ$, such that the module $V_{\om_i}(z_i)\ox V_{\om_j}(z_j)$
is reducible.
In the special case when both $\om_i$ and $\om_j$ are usual Young
diagrams, this description was recently given by A.\,Molev in [M].

The {\it only if\/} part of Theorem~4.8
is an easy and well known fact.
Our proof of the {\it if\/} part is based on the
following general observation. Let $\A$ be any
unital associative algebra over $\CC$,
and let $W$ be any finite-dimensional $\A\]$-module.
Denote by $\rho$ the corresponding homomorphism $\A\to\End(W)$.
Suppose that the subalgebra
$$
\rho\ts(\A)\ox\End(W)\subset\End(W)\ox\End(W)=\End(W\ox W)
$$
contains the flip map $x\ox y\mapsto y\ox x$ in $W\ox W$.
Then $\rho\ts(\A)=\End(W)$ and
the $\A\]$-module $W$ is irreducible, see Lemma~4.10.
We employ this observation
when $\A=\YN$, and $W$ is the tensor product (0.2). Moreover,
to prove the irreducibility of (0.2) under
the condition that all the operators $I_{ij}$ with
$i<j$ are invertible, it suffices to show that
$I_{ii}$ is proportional to the flip map in $V_{\om_i}\ox V_{\om_i}$
for each index~$i$. For the details of this argument see Proposition~4.14.
We took the idea of this argument from the recent physical
literature [\ts KMT,\ts MT\ts]\ts. This argument extends to
the quantum affine algebra
$\text{U}_{\ns q}(\widehat{\frak{\ts gl\ts }}_{\ns N})$
in a straightforward way. We work with the Yangian $\YN$
in this article only
to simplify the exposition. Moreover, this argument
extends to other quantum affine algebras due to the existence
of the universal $R\]$-matrix\ts; see \text{[\ts FM\ts,\ts K\ts]}
for the precise formulation of the problem and known results in
the general case.

In this article,
we use the explicit realization of the elementary $\YN\]$-module $V_\om(z)$,
introduced in [C2]. Denote $\la_1-\mu_1+\la_2-\mu_2+\ldots=n$.
Let $n>0$.
Consider the action of the symmetric group $S_n$ in
$(\CC^N)^{\ox n}$ by permutations of the $n$ tensor factors.
Denote by $\pi_n$
the corresponding homomorphism $\CS_n\to\End\bigl((\CC^N)^{\ox n}\bigr)$.
Then $V_\om$ is the subspace in $(\CC^N)^{\ox n}\]$
defined as the image of the operator $\pi_n(F_\om)$ for a certain element
$F_\om$ of the group ring $\CS_n$. This element is constructed
by a limiting process, called the fusion procedure [C1].
Note that in [C1] most of the results were given without proofs.
In Section~2 we give proofs of all properties of the element $F_\om$
which we need in this article. These proofs are new, cf.~[N].

This choice of the realization of the $\YN\]$-module $V_\om(z)$
allows us to write down explicitly,
for any $z'\in\CC$ such that $z-z'\notin\ZZ$,
an invertible intertwining operator
$$
R(z,\ns z'):\ts
V_\om(z)\ox V_\om(z')\longrightarrow
V_\om(z)\oxp V_\om(z')\,.
$$
This operator comes with a natural normalization.
In partucular, it is a rational function in $z-z'$ with values
in $\End(V_\om\ox V_\om)$. Denote by $P_\om$ the flip map in
$V_\om\ox V_\om$. We prove that
the leading coefficient of the Laurent expansion in $u$
near the origin $u=0$
of the function $R(z+\ns u\ts,\ns z)$ equals
$c(\om)\ts P_\om\,u^{-d(\om)}\]$ for a certain non-zero rational
number $c(\om)$ and a positive integer $d(\om)$,
see Proposition~4.7. The proof of this proposition is
based on the results of Section~3.
We call $d(\om)$ the Durfee rank of the skew
Young diagram $\om$. When $\om$ is a usual Young diagram,
$d(\om)$ is the number of boxes on its main diagonal.
In Section~1 we study the remarkable
combinatorial properties of the numbers $d(\om)$
for all skew Young diagrams $\om$.


\Section{Durfee rank of a skew Young diagram}
\nt
Let $\la=(\la_1,\la_2\ts,\,\ldots)$ and $\mu=(\mu_1,\mu_2\ts,\,\ldots)$
be any two partitions. As usual, here the parts are arranged in the
non-increasing order\ts: we have $\la_1\ge\la_2\ge\ldots\ge0$ and
$\mu_1\ge\mu_2\ge\ldots\ge0$. We will always
suppose that $\la_i\ge\mu_i$ for any
$i=1,2,\ts\ldots\ts$.
Consider the {\it skew Young diagram} $\om=\la\ts/\mu$.
It can be defined as the set of pairs
$$
\{\,(i,\ns j)\in\ZZ^2\ |\ i\ge1,\ \la_i\ge j>\mu_i\,\}
$$
When $\mu=(0,\ns0\ts,\ts\ldots)$,
this is the Young diagram of the partition
$\la$. We will employ the standard graphic representation of Young
diagrams on the plane $\RR^2$ with the matrix style coordinates $(x,\ns y)$.
Here the first coordinate $x$ increases from top to bottom, while the second
coordinate $y$ increases from left to right. The element $(i,\ns j)\in\om$
is represented by the unit box with the bottom right corner
at the point $(i,\ns j)\in\RR^2$. For example,
below we represent the skew Young diagram $\om$ for the partitions
$\la=(9,\ns9,\ns9,\ns7,\ns7,\ns3,\ns3,\ns3,\ns3,\ns0,\ns0,\ts\ldots)$ and
$\mu=(5,\ns5,\ns3,\ns3,\ns3,\ns3,\ns2,\ns0,\ns0,\ts\ldots)$.
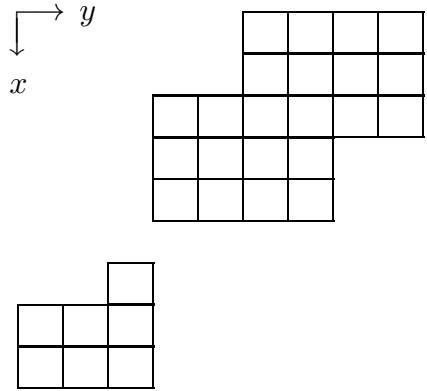
\bigskip
\vbox{
$$
\kern23pt\longrightarrow\,y\mv
$$
\vglue-26.182pt
$$
\kern-1pt\vert\mv
$$
\vglue-28pt
$$
\bigr\downarrow\mv\kern0.991pt
$$
\vglue-16pt
$$
x\mv
$$
\vglue-46.5pt
$$
\phantom{\bx}
\phantom{\bx}
\phantom{\bx}
\phantom{\bx}
\phantom{\bx}
{\bx}
{\bx}
{\bx}
{\bx}
$$
\vglue-17.8pt
$$
\phantom{\bx}
\phantom{\bx}
\phantom{\bx}
\phantom{\bx}
\phantom{\bx}
{\bx}
{\bx}
{\bx}
{\bx}
$$
\vglue-17.8pt
$$
\phantom{\bx}
\phantom{\bx}
\phantom{\bx}
{\bx}
{\bx}
{\bx}
{\bx}
{\bx}
{\bx}
$$
\vglue-17.8pt
$$
\phantom{\bx}
\phantom{\bx}
\phantom{\bx}
{\bx}
{\bx}
{\bx}
{\bx}
\phantom{\bx}
\phantom{\bx}
$$
\vglue-17.8pt
$$
\phantom{\bx}
\phantom{\bx}
\phantom{\bx}
{\bx}
{\bx}
{\bx}
{\bx}
\phantom{\bx}
\phantom{\bx}
$$
\vglue-17.8pt
$$
\phantom{\bx}
\phantom{\bx}
\phantom{\bx}
\phantom{\bx}
\phantom{\bx}
\phantom{\bx}
\phantom{\bx}
\phantom{\bx}
\phantom{\bx}
$$
\vglue-17.8pt
$$
\phantom{\bx}
\phantom{\bx}
{\bx}
\phantom{\bx}
\phantom{\bx}
\phantom{\bx}
\phantom{\bx}
\phantom{\bx}
\phantom{\bx}
$$
\vglue-17.8pt
$$
{\bx}
{\bx}
{\bx}
\phantom{\bx}
\phantom{\bx}
\phantom{\bx}
\phantom{\bx}
\phantom{\bx}
\phantom{\bx}
$$
\vglue-17.8pt
$$
{\bx}
{\bx}
{\bx}
\phantom{\bx}
\phantom{\bx}
\phantom{\bx}
\phantom{\bx}
\phantom{\bx}
\phantom{\bx}
$$
}
\bigskip
\nt
We will identify the elements of the set $\om$ with the corresponding unit
boxes~on~$\RR^2$.

Let us call an element $(i,\ns j)\in\om$ {\it left-convex\/}
if $(i\!-\!\ns1,\ns j)\notin\om$ and $(i,\ns j\!-\!\ns1)\notin\om$.
An element $(i,\ns j)\in\om$ will be called {\it left-concave\/} if both
$(i\!-\!\ns1,\ns j)\in\om$ and $(i,\ns j\!-\!\ns1)\in\om$,
but $(i\!-\!\ns1,\ns j\!-\!\ns1)\notin\om$.
Each of these two definitions has a natural interpretation in terms of
boxes of the Young diagram $\om$. If $\la\neq\mu=(0,\ns0\ts,\ts\ldots)$ then
there is only one left-convex element $(1,\ns1)\in\om$
and no left-concave elements. Now consider those diagonals of $\om$,
which contain a left-convex box. Let $d$ be the total number of boxes
on those diagonals. Let $\db$
be the total number of boxes on the diagonals of $\om$
which contain a left-concave box. The difference $d-\db$ will be called
the {\it Durfee rank} of the skew diagram $\om$ and denoted by $d(\om)$.
If $\mu=(0,\ns0\ts,\ts\ldots)$ then $d(\om)=d$ is the number of
boxes on the main diagonal
$
\{\,(i,\ns i)\in\ZZ^2\ |\ \la_i\ge i\ge1\,\}
$
of the Young diagram of $\la$, and our definition coincides with the usual
one \text{[\,A,\,Section\,2.3\,]}.
Below for our exemplary skew Young diagram $\om$ we mark
by the $+$ or $-$
signs the diagonals, which
contain a left-convex or a left-concave box respectively.
\bigskip\medskip
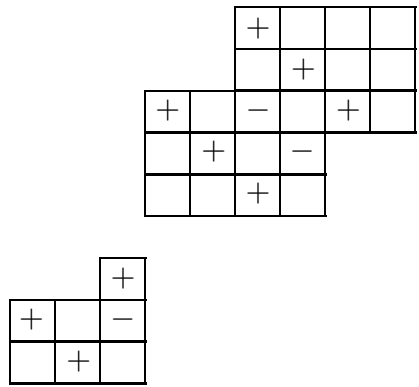
\vbox{
$$
\phantom{\bx}
\phantom{\bx}
\phantom{\bx}
\phantom{\bx}
\phantom{\bx}
{\bx}
{\bx}
{\bx}
{\bx}
$$
\vglue-17.8pt
$$
\phantom{\bx}
\phantom{\bx}
\phantom{\bx}
\phantom{\bx}
\phantom{\bx}
{\bx}
{\bx}
{\bx}
{\bx}
$$
\vglue-17.8pt
$$
\phantom{\bx}
\phantom{\bx}
\phantom{\bx}
{\bx}
{\bx}
{\bx}
{\bx}
{\bx}
{\bx}
$$
\vglue-17.8pt
$$
\phantom{\bx}
\phantom{\bx}
\phantom{\bx}
{\bx}
{\bx}
{\bx}
{\bx}
\phantom{\bx}
\phantom{\bx}
$$
\vglue-17.8pt
$$
\phantom{\bx}
\phantom{\bx}
\phantom{\bx}
{\bx}
{\bx}
{\bx}
{\bx}
\phantom{\bx}
\phantom{\bx}
$$
\vglue-17.8pt
$$
\phantom{\bx}
\phantom{\bx}
\phantom{\bx}
\phantom{\bx}
\phantom{\bx}
\phantom{\bx}
\phantom{\bx}
\phantom{\bx}
\phantom{\bx}
$$
\vglue-17.8pt
$$
\phantom{\bx}
\phantom{\bx}
{\bx}
\phantom{\bx}
\phantom{\bx}
\phantom{\bx}
\phantom{\bx}
\phantom{\bx}
\phantom{\bx}
$$
\vglue-17.8pt
$$
{\bx}
{\bx}
{\bx}
\phantom{\bx}
\phantom{\bx}
\phantom{\bx}
\phantom{\bx}
\phantom{\bx}
\phantom{\bx}
$$
\vglue-17.8pt
$$
{\bx}
{\bx}
{\bx}
\phantom{\bx}
\phantom{\bx}
\phantom{\bx}
\phantom{\bx}
\phantom{\bx}
\phantom{\bx}
$$
\vglue-137.5pt
$$
\kern29pt+
$$
\vglue-18pt
$$
\kern57.5pt+
$$
\vglue-18pt
$$
\kern29pt+\kern18.5pt-\kern18.5pt+
$$
\vglue-18pt
$$
\kern29pt+\kern20pt-
$$
\vglue-18pt
$$
\kern29pt+
$$
\vglue-4pt
$$
+\kern56pt
$$
\vglue-18pt
$$
+\kern21pt-\kern85pt
$$
\vglue-18pt
$$
+\kern84pt
$$
}
\bigskip\noindent
In this example the Durfee rank $d(\om)$ equals $9-3=6$.
Our first observation is

\proclaim{Proposition 1.1}
If the diagram $\om$ is not empty then $d(\om)>0$.
\endproclaim

\demo{Proof}
For each $i=1,\ns2,\ts\ldots$
let $d_i$ be the number of boxes on that diagonal of $\om$
which contains the leftmost box of the $i\]$-th row.
If the $i\]$-th row of $\om$ is empty, we set $d_i=0$.
Let $\db_i$ be the number of boxes on the next diagonal down. Then
$$
d(\om)=d_1-\db_1+d_2-\db_2+\ldots\,\ts.
\Tag{1.1}
$$
By the definition of a skew Young diagram
we have $d_i-\db_i\ge0$ for every index $i=1,\ns2,\ts\ldots\ts$.
But if $i$ is the index of the last non-empty row of $\om$,
then $d_i=1$ and $\db_i=0$
\enddemos

\noindent
Denote by $\ell(\om)$ the number of non-empty rows~in~$\om$.
Here is a formula for $d(\om)$.

\proclaim\nofrills{Proposition 1.2}{\bf\,:\ }
$
d(\om)=\ell(\om)-
\#\ts\{\,(i,\ns j)\ |\ \la_j\ns-\ns j=\mu_i\ns-\ns i\ts,\ i<j\ts\}\,.
$
\endproclaim

\demo{Proof}
We will keep to the notation from the proof of Proposition 1.1. Further,
for each $i=1,\ns2,\ts\ldots$
put $e_i=0$ or $e_i=1$ depending on whether the $i\]$-th row of $\om$
is empty or not. By the definition of a skew Young diagram
we have the equalities
$$
\align
d_i\ts=\ts{}&\#\ts
\{\,j\ |\ \la_j\ns-\ns j\ge\mu_i\ns-\ns i\ns+\ns1\ge\mu_j\ns-\ns j\ns+\ns1\,\}
\\
\ts=\ts{}&\#\ts
\{\,j\ |\ \la_j\ns-\ns j\ge\mu_i\ns-\ns i\ns+\ns1\ts,\ i\le j\ts\}
\\
\ts=\ts{}&\#\ts
\{\,j\ |\ \la_j\ns-\ns j\ge\mu_i\ns-\ns i\ns+\ns1\ts,\ i<j\ts\}+e_i\,,
\\
\db_i\ts=\ts{}&\#\ts
\{\,j\ |\ \la_j\ns-\ns j\ge\mu_i\ns-\ns i\ge\mu_j\ns-\ns j\ns+\ns1\,\}
^{\phantom{|}}
\\
\ts=\ts{}&\#\ts
\{\,j\ |\ \la_j\ns-\ns j\ge\mu_i\ns-\ns i\ts,\ i<j\ts\}\,.
\endalign
$$
Therefore,
$$
d_i-\db_i\ts=
\ts e_i-\#\ts\{\,j\ |\ \la_j\ns-\ns j=\mu_i\ns-\ns i\ts,\ i<j\ts\}\,.
$$
Hence we obtain from \(1.1) that
$$
d(\om)=(e_1+e_2+\ldots)-
\#\ts\{\,(i,\ns j)\ |\ \la_j\ns-\ns j=\mu_i\ns-\ns i\ts,\ i<j\ts\}\,.
$$
By observing that $e_1+e_2+\ldots=\ell(\om)$ we now complete the proof
\enddemos

\noindent
Let $\lap=(\lap_1,\lap_2\ts,\,\ldots)$ and $\mup=(\mup_1,\mup_2\ts,\,\ldots)$
be the partitions conjugate to $\la$ and $\mu$ respectively.
Consider the corresponding skew Young diagram $\omp=\lap\!/\mup$.
It is obtained from the diagram $\om\subset\ZZ^2$ by transposition
$(i,\ns j)\mapsto(j,\ns i)$. Evidently, $d(\omp)=d(\om)$.
Therefore,
Proposition 1.2 can be reformulated as

\proclaim\nofrills{Proposition 1.3}{\bf\,:\ }
$
d(\om)=\ell(\omp)-
\#\ts\{\,(i,\ns j)\ |\ \lap_j\ns-\ns j=\mup_i\ns-\ns i\ts,\ i<j\ts\}\,.
$
\endproclaim

\noindent
By our definition, the Durfee rank $d(\om)$ depends only on the shape of
the diagram $\om$. That is, if a skew Young diagram $\omf$
is obtained from $\om\subset\ZZ^2$ by translation
$(i,\ns j)\mapsto(i\ns+\ns k,\ns j\ns+\ns l)$ for a certain
$(k,\ns l)\in\ZZ^2$, then $d(\omf)=d(\om)$. Furthermore, suppose that
a skew Young diagram $\omn$ is obtained from $\om$ via the map
$(i,\ns j)\mapsto(k\ns-\ns i\ns+\ns1,\ns l\ns-\ns j\ns+\ns1)$.
The boxes of $\omn$ on the plane $\RR^2$ are obtained from those
of $\om$ via rotation by $180^{\!\!\crc\!\!\!}$
around the point $(\frac{k}2,\frac{l}2)$. The main result of this section is

\proclaim\nofrills{Theorem 1.4}{\bf\,:\ }
$d(\omn)=d(\om)$.
\endproclaim

\demo{Proof}
We have the inequalities $\la_1\ns-\ns1>\la_2\ns-\ns2>\ldots$
and $\mu_1\ns-\ns1>\mu_2\ns-\ns2>\ldots$\,. So under the condition
$\la_j\ns-\ns j=\mu_i\ns-\ns i\ts$ in Proposition 1.2
we have the equivalencies
$$
\alignat4
&i<j
&\ \ \Leftrightarrow
&\ \ \la_i\ns-\ns i>\la_j\ns-\ns j
&\ \ \Leftrightarrow
&\ \ \la_i\ns-\ns i>\mu_i\ns-\ns i
&\ \ \Leftrightarrow
&\ \ \la_i>\mu_i\,,
\\
&i<j
&\ \ \Leftrightarrow
&\ \ \mu_i\ns-\ns i>\mu_j\ns-\ns j
&\ \ \Leftrightarrow
&\ \ \la_j\ns-\ns j>\mu_j\ns-\ns j
&\ \ \Leftrightarrow
&\ \ \la_j>\mu_j\,.
\endalignat
$$
Denote by $\Cal L(\om)$ the set of differences
$\mu_i\ns-\ns i$ such that the $i\]$-th row of $\om$ is not empty.
The set of differences
$\la_j\ns -\ns j$ such that the $j\]$-th row of $\om$ is not empty,
will be denoted by $\R(\om)$. By the above equivalencies,
Proposition 1.2 can be restated~as
$$
d(\om)=\ell(\om)-\#\ts\bigl(\L(\om)\cap\R(\om)\bigr)\,.
\Tag{1.2}
$$
We have $\ell(\omn)=\ell(\om)$. We also have
$\#\ts\bigl(\L(\omn)\cap\R(\omn)\bigr)=\#\ts\bigl(\L(\om)\cap\R(\om)\bigr)$
since
$$
\L(\omn)=l\ns-\ns k\ns-\ns\R(\om)-\!1\ts,
\quad
\R(\omn)=l\ns-\ns k\ns-\ns\L(\om)-\!1\ts,
$$
by the defintion of $\omn$.
Thus we obtain Corollary 1.4 from the formula \(1.2)
\enddemos

\noindent
Let us call an element $(i,\ns j)\in\om$ {\it right-convex\/}
if $(i\!+\!\ns1,\ns j)\notin\om$ and $(i,\ns j\!+\!\ns1)\notin\om$.
An element $(i,\ns j)\in\om$ is called {\it right-concave\/} if both
$(i\!+\!\ns1,\ns j)\in\om$ and $(i,\ns j\!+\!\ns1)\in\om$,
but $(i\!+\!\ns1,\ns j\!+\!\ns1)\notin\om$.
Each of these two definitions
has a natural interpretation in terms of
boxes of the Young diagram $\om$.
We could define the Durfee rank of $\om$ using the
right-convex and right-concave boxes instead of the
left-convex and left-concave, respectively.
Theorem 1.4 ensures that we would then get the same number $d(\om)$.
Note that although the
Young diagram of a non-empty partition always
has only one left-convex box and no left-concave boxes, it may have
several right-convex and right-concave boxes.
\text{So our definition of the Durfee rank of $\om$ is quite natural.}


\Section{Fusion procedure for a skew Young diagram}
\nt
Take a non-empty skew Young diagram $\om$.
Let $n$ be the total number of boxes in $\om$.
In this section we will introduce a certain element $F_\om$ of
the symmetric group ring $\CS_n$. Under some extra conditions on
the diagram $\om$, this element will be used in Section 4
to determine a family of irredicible modules $V_\om(z)$
over the Yangian $\YN$. Here the parameter $z$ is ranging over
the complex field $\CC$. For the skew Young diagram $\omf$
obtained from $\om$ by translation
$(i,\ns j)\mapsto(i\ns+\ns k,\ns j\ns+\ns l)$ for a certain
$(k,\ns l)\in\ZZ^2\]$, we will have $V_{\omf}(z)=V_\om(z-k+l)$.

Consider the {\it column tableau\/} of shape $\om$.
It is obtained by filling
the boxes of the diagram $\om$ with the numbers $1\lcd n$
consecutively by columns from left to right, downwards in every column.
We will denote this tableau by $\Om$.
Further, for each $p=1\lcd n$ put $c_p=j-i$ if the box
$(i,j)\in\om$ is filled with the number $p$ in the column tableau $\Om$.
The difference $j-i$ is called the {\it content\/} of the box $(i,j)$
of the diagram $\om$. Our choice of the
tableau $\Om$ provides an ordering of the collection of all contents of $\om$.
Below on the left for the partitions
$\la=(5,\ns3,\ns3,\ns3,\ns3,\ns0,\ns0,\ts\ldots)$ and
$\mu=(3,\ns3,\ns3,\ns2,\ns0,\ns0,\ts\ldots)$
we show the column tableau $\Om$. On the right we
indicate the contents of all boxes of the diagram $\om$.
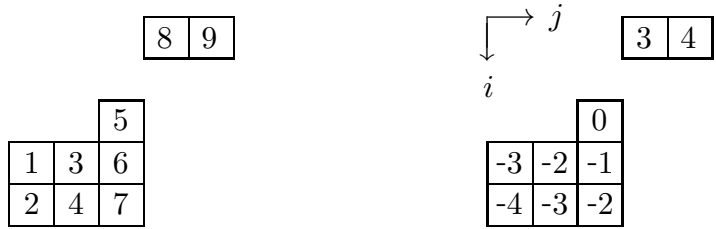
\bigskip
\vbox{
$$
\kern22.5pt\longrightarrow\,j\mw
$$
\vglue-26.182pt
$$
\kern-1pt\vert\mw
$$
\vglue-28pt
$$
\bigr\downarrow\mw\kern0.991pt
$$
\vglue-16pt
$$
i\mw
$$
\vglue-46.5pt
$$
\phantom{\bx}
\phantom{\bx}
\phantom{\bx}
{\bx}
{\bx}
\kern80pt
\phantom{\bx}
\phantom{\bx}
\phantom{\bx}
{\bx}
{\bx}
$$
\vglue-17.8pt
$$
\phantom{\bx}
\phantom{\bx}
\phantom{\bx}
\phantom{\bx}
\phantom{\bx}
\kern80pt
\phantom{\bx}
\phantom{\bx}
\phantom{\bx}
\phantom{\bx}
\phantom{\bx}
$$
\vglue-17.8pt
$$
\phantom{\bx}
\phantom{\bx}
{\bx}
\phantom{\bx}
\phantom{\bx}
\kern80pt
\phantom{\bx}
\phantom{\bx}
{\bx}
\phantom{\bx}
\phantom{\bx}
$$
\vglue-17.8pt
$$
{\bx}
{\bx}
{\bx}
\phantom{\bx}
\phantom{\bx}
\kern80pt
{\bx}
{\bx}
{\bx}
\phantom{\bx}
\phantom{\bx}
$$
\vglue-17.8pt
$$
{\bx}
{\bx}
{\bx}
\phantom{\bx}
\phantom{\bx}
\kern80pt
{\bx}
{\bx}
{\bx}
\phantom{\bx}
\phantom{\bx}
$$
\vglue-84pt
$$
\kern42pt8\kern9pt9\kern90pt\kern42pt3\kern9pt4
$$
\vglue-5pt
$$
5\kern146pt0
$$
\vglue-17.5pt
$$
1\kern9pt3\kern9pt6\kern116pt\mi3\kern6pt\mi2\kern6pt\mi1\kern26pt
$$
\vglue-18pt
$$
2\kern9pt4\kern9pt7\kern116pt\mi4\kern6pt\mi3\kern6pt\mi2\kern26pt
$$
}
\bigskip\noindent
Here $n=9$ and the sequence of contents $(c_1\lcd c_{\ts9})$ is
$(\ts\mi3,\ns\mi4,\ns\mi2,\ns\mi3,\ns0,\ns\mi1,\ns\mi2,\ns3,\ns4\ts)$.

For any two distinct numbers $p,\ns q\in\{1\lcd n\}$
let $(p\ts q)$ be the
transposition in the symmetric group $S_n$.
Consider the rational functions of two complex variables $u\ts,\ns v$
with values in the group ring $\CS_n$
$$
f_{pq}(u,\ns v)=1-\frac{(p\ts q)}{u-v}\,.
\Tag{2.0}
$$
As a direct calculation shows, these rational functions satisfy the equations
$$
f_{pq}(u,\ns v)\ts f_{pr}(u,\ns w)\ts f_{qr}(v,\ns w)=
f_{qr}(v,\ns w)\ts f_{pr}(u,\ns w)\ts f_{pq}(u,\ns v)
\Tag{2.1}
$$
for all pairwise distinct indices $p,\ns q,\ns r$. Evidently,
for all pairwise distinct $p,\ns q,\ns s,\ns t$
$$
f_{pq}(u,\ns v)\ts f_{st}(z,\ns w)=f_{st}(z,\ns w)\ts f_{pq}(u,\ns v)\,.
\Tag{2.2}
$$
For all distinct $p,q$ we also have
$$
f_{pq}(u,\ns v)\,f_{qp}(v,\ns u)=1-\frac1{(u-v)^2}\,.
\Tag{2.25}
$$
Our construction of the element $F_\om\in\CS_n$ will be based on the following
simple observation. Consider the
rational function of $u\ts,\ns v\ts,\ns w$ defined by
the product at either side of \(2.1). The factor
$f_{pr}(u,\ns w)$ in \(2.1) has a pole at $u=w$. However, we have

\proclaim{Lemma 2.1}
The restriction of the rational function \(2.1) to the set of
$(u\ts,\ns v\ts,\ns w)$ such that $u=v\pm1$, is regular at $u=w$.
\endproclaim

\demo{Proof}
Under the condition $u=v\pm1$ the product on the left hand side of \(2.1)
can be written as
$$
\bigl(\ts1\mp(p\ts q)\bigr)\cdot\bigg(1-\frac{(p\ts r)+(q\ts r)}{v-w}\bigg).
$$
The latter rational function of $v,\ns w$ is manifestly regular at
$w=v\pm1$
\enddemos

\noindent
Now introduce $n$ complex variables $z_1\lcd z_n$.
Equip the set of all pairs $(p,\ns q)$
where $1\le p<q\le n$, with the lexicographical ordering.
Take the ordered product
$$
\prod_{(p,q)}^{\longrightarrow}\ f_{pq}(c_p+z_p\ts,c_q+z_q)
\Tag{2.3}
$$
over this set. Consider \(2.3)
as a rational function of the variables $z_1\lcd z_n$ with values
in $\CS_n$. Denote this rational function by $F_\om(z_1\lcd z_n)$.
Let $\Cal D_\om$ be the vector subspace in $\CC^n$ consisting of all tuples
$(z_1\lcd z_n)$ such that $z_p=z_q$ whenever the numbers $p$ and $q$
appear in the same column of the tableau $\Om$.
The origin $\or\in\CC^n$ belongs to $\D_\om$.
The following statement goes back to [C1].

\proclaim{Proposition 2.2} The restriction of the rational
function $F_\om(z_1\lcd z_n)$ to the subspace $\D_\om\subset\CC^n$
is regular~at the point $\or$.
\endproclaim

\demo{Proof}
We shall provide an expression for the restriction of
$F_\om(z_1\lcd z_n)$ to $\D_\om$ which is manifestly regular at $\or$.
The factor $f_{pq}(c_p+z_p\ts,c_q+z_q)$ has a pole at $z_p=z_q$ if and
only if the numbers $p$ and $q$ stand on the same diagonal of the
tableau $\Om$. We shall then call the pair $(p,\ns q)$ {\it singular.}

Let any singular pair $(p,\ns q)$ be fixed.
We will write $(p,\ns q)\prec(s,\ns t)$ if the pair
$(p,\ns q)$ precedes $(s,\ns t)$ in the lexicographical ordering, that is,
if $p<s$ or if $p=s$,$q<t$. In this ordering,
the pair immediately before $(p,\ns q)$ is $(p,\ns q-1)$.
Moreover, the number $q-1$ stands just above $q$ in the
tableau $\Om$. Therefore,
$z_{q-\ns1}=z_q$ on $\D_\om$.
We also have $(p,\ns q)\prec(q\!-\!1,\ns q)$, since $p<q$.
By \(2.1) and \(2.2) the product over the pairs $(s,\ns t)$
$$
\prod_{(p,q)\prec(s,t)}^{\longrightarrow}\ f_{st}(c_s+z_s\ts,c_t+z_t)
$$
is divisible on the left
by the factor $f_{q-\ns1,q}(c_{q-\ns1}\ns+\ns z_{q-\ns1}\ts,c_q\ns+\ns z_q)$.
Restriction of this factor to $z_{q-\ns1}=z_q$
is just $1\ns-(q-1,\ns q)$ since $c_{q-\ns1}=c_q\ns+\ns1$.

The element $\bigl(1-(q-1,\ns q)\bigr)/2\in\CS_n$ is an idempotent.
Now for every singular pair $(p,\ns q)$ let us replace the two adjacent
factors in \(2.3)
$$
f_{p,q-\ns1}(c_p+z_p\ts,c_{q-\ns1}+z_{q-\ns1})\ts
f_{pq}(c_p+z_p\ts,c_q+z_q)
$$
\vskip-10pt\noindent
by
$$
f_{p,q-\ns1}(c_p+z_p\ts,c_{q-\ns1}+z_{q-\ns1})\ts
f_{pq}(c_p+z_p\ts,c_q+z_q)\ts
f_{q-\ns1,q}(c_{q-\ns1}+z_{q-\ns1},c_q+z_q)/2\,.
$$
This does not affect the values of the restriction of \(2.3) to $\D_\om$.
But the restriction to $z_{q-\ns1}=z_q$ of the replacement product
is regular at $z_p=z_q$ by Lemma 2.1
\enddemos

\noindent
Due to Proposition 2.2,
an element $F_\om\in\CS_n$ can be now defined as the value at
the point $\or$
of the restriction to $\D_\om$ of the function $F_\om(z_1\lcd z_n)$.
Note that for $n=1$ we have $F_\om=1$. For any $n\ge1$ we have
the following fact.

\proclaim{Proposition 2.3}
The coefficient of $F_\om\in\CS_n$ at the unit element of
$S_n$ is $1$.
\endproclaim

\demo{Proof}
For each $r=1\lcd n-1$ let $g_r\in S_n$ be the adjacent transposition
$(r,r\ns+\ns 1)$. Let $g_0\in S_n$ be the element of the maximal length.
Let us multiply the ordered product \(2.3) by the element $g_0$ on the right.
Using the reduced decomposition
$$
g_0\,=\prod_{(p,q)}^{\longrightarrow}g_{q-p}\ ,
\Tag{2.31}
$$
\vskip-7.5pt\noindent
we get the product
$$
\prod_{(p,q)}^{\longrightarrow}\
\left(g_{q-p}-\frac1{c_p\ns+\ns z_p\ns-\ns c_q\ns-\ns z_q}\right).
$$
It expands as a sum of the elements $g\in S_n$ with the coefficients
from the field of rational functions of $z_1\lcd z_n$ valued in $\CC$.
Since the decomposition \(2.31) is reduced, the coefficient at $g_0$
is $1$. By the definition of $F_\om$
this implies that
the coefficient of $F_\om\ts g_0\in\CS_n$ at
$g_0\in S_n$ is also $1$
\enddemos

\noindent
In particular, Proposition 2.3 shows that $F_\om\neq0$ for any
non-empty diagram $\om$.
Let us now denote by $\al$ the involutive antiautomorphism
of the group ring $\CS_n$ defined by $\al(g)=g^{-1}$ for every $g\in S_n$.

\proclaim{Proposition 2.4}
The element $F_\om\in\CS_n$ is $\al\]$-invariant.
\endproclaim

\demo{Proof}
Any element of the group ring $\CS_n$ of the form
$f_{pq}(u,\ns v)$ is $\al$\ns-invariant\,.
Applying the antiautomorphism $\al$ to an element of the form
\(2.3) just reverses the ordering of
the factors corresponding to the
pairs $(p,\ns q)$. However, the initial ordering can be then restored by
using the relations \(2.1) and \(2.2). Therefore,
any value of the function $F_\om(z_1\lcd z_n)$ is $\al$\ns-invariant.
So is the element $F_\om$
\enddemos

\noindent
The next property of the element $F_\om\in\CS_n$
also easily follows from its definition.

\proclaim{Theorem 2.5}
Suppose the numbers $r$ and $r\ns+\ns1$ stand in the
same column of $\Om$. Then the element $F_\om\in\CS_n$
is divisible by $f_{r,r+1}(\ts c_r,\ns c_{r+1})=1\ns-\ns(r,r\!+\!1)$
on the left and right.
\endproclaim

\demo{Proof}
By Proposition 2.4 the divisibility of $F_\om$ by
the element $1\!-\ns(r,r\!+\!1)$ on the left is equivalent
to the divisibility by the same element on the right.
Let us prove the divisibility on the left.
Using the definition \(2.3) along with the
relations \(2.1) \ns and \ns\(2.2), we can write
$$
F_\om(z_1\lcd z_n)=
f_{r,r+1}(\ts c_r\ns+\ns z_r,\ns c_{r+1}\ns+\ns z_{r+1}\ts)\,
F(z_1\lcd z_n)
\Tag{2.33}
$$
for some rational function $F(z_1\lcd z_n)$ valued in $\CS_n$.
By Proposition 2.2, the restriction of the function
$F_\om(z_1\lcd z_n)$ to $\D_\om$
is regular at $\or$. But the restriction to $\D_\om$ of the factor
$f_{r,r+1}(\ts c_r\ns+\ns z_r,\ns c_{r+1}\ns+\ns z_{r+1}\ts)$ in \(2.33)
is $f_{r,r+1}(\ts c_r,\ns c_{r+1})$.
In particular, the value at $\or$
of the restriction of $F_\om(z_1\lcd z_n)$
to $\D_\om$ is divisible on the left by
$f_{r,r+1}(\ts c_r,\ns c_{r+1})=1\ns-\ns(r,r\!+\!1)$
\enddemos

\noindent
Let $l=1\lcd n-1$. Denote by $\si_l$ the embedding
$\CS_{n-l}\to\CS_n$ defined by
$\si_l:(p\ts q)\mapsto(p\ns+\ns l,\ns q\ns+\ns l)$
for all distinct $p,q=1\lcd n-l$.
Define a skew Young diagram $\ph$ as the shape of the tableau obtained
from $\Om$ by removing each of the numbers $1\lcd l$. Consider
the corresponding element $F_\ph\in\CS_{n-l}$.

\proclaim{Proposition 2.6}
The element $F_\om$ is divisible by $\si_l(F_\ph)$
on the left and right.
\endproclaim

\demo{Proof}
Due to Proposition 2.4 the divisibility of $F_\om$ by
the element $\si_l(F_\ph)$ on the left is equivalent
to the divisibility by the same element on the right.
We will actually prove the divisibility on the right.
By the definition \(2.3) we have
$$
F_\om(z_1\lcd z_n)=
\prod_{(p,q)}^{\longrightarrow}\ f_{pq}(c_p+z_p\ts,c_q+z_q)
\,\cdot\,
\si_l\bigl(F_\ph(z_{l+1}\lcd z_n)\bigr)
\Tag{2.333}
$$
where the pairs $(p,\ns q)$ are ordered lexicigraphically, and
precede the pair $(l\ns+\ns1,\ns l\ns+\ns2)$.
According to our proof of Proposition 2.2, the product
over the pairs $(p,\ns q)$ in \(2.333) maybe replaced
by a certain rational function of $z_1\lcd z_n$ with the restriction
to $\D_\om$ regular at $\or$, without affecting the values of the
restriction. Now the divisibility of $F_\om$ by $\si_l(F_\ph)$
on the right follows from the decomposition \(2.333)
and the definition of the element $F_\ph\in\CS_{n-l}$
\enddemos

\noindent
Let $m=1\lcd n-1$.
Regard the group ring $\CS_m$ as a subalgebra in $\CS_n$
using the standard embedding $S_m\to S_n$.
Define a skew Young diagram $\ps$ as the shape of the tableau obtained
from $\Om$ by removing each of the numbers $m\ns+\ns1\lcd n$. Consider
the corresponding element $F_\ps\in\CS_m$.
The following property of the element $F_\om$ is obtained by combining the
arguments of Propositions 2.2 and 2.6.

\proclaim{Proposition 2.7}
The element $F_\om\in\CS_n$ is divisible by $F_\ps$ on the left
and right.
\endproclaim

\demo{Proof}
Let us introduce another expression for restriction of
$F_\om(z_1\lcd z_n)$ to $\D_\om$ which is again manifestly regular at $\or$.
Reorder the pairs $(p,\ns q)$ in the product \(2.3) as follows;
this reordering will not affect the values of the
product due to \(2.1) and \(2.2).
The pair $(s,\ns t)$ will precede the pair $(p,\ns q)$
if $t<q$ or if $t=q$,$s<p$. Let us now write
$(s,\ns t)\prec(p,\ns q)$ referring to this new ordering.

Let any singular pair $(p,\ns q)$ be fixed.
In our new ordering,
the pair immediately after $(p,\ns q)$ is $(p\!+\!1,\ns q)$.
Moreover, the number $p\!+\!1$ stands just below $p$ in the
tableau $\Om$. Therefore,
$z_p=z_{p+\ns1}$ on $\D_\om$.
We also have $(p,\ns p\!+\!1)\prec(p,\ns q)$, since $p<q$.
By \(2.1) and \(2.2) the product
$$
\prod_{(s,t)\prec(p,q)}^{\longrightarrow}\ f_{st}(c_s+z_s\ts,c_t+z_t)
$$
over the pairs $(s,\ns t)$
is divisible on the right by the factor
$f_{p,p+\ns1}(c_p\ns+\ns z_p\ts,c_{p+\ns1}\ns+\ns z_{p+\ns1})$.
Restriction of this factor to $z_p=z_{p+\ns1}$
is just $1\ns-(p,\ns p+1)$ since $c_{p+\ns1}=c_p\ns-\ns1$.

The element $\bigl(1-(p,\ns p+1)\bigr)/2\in\CS_n$ is an idempotent.
Now for every singular pair $(p,\ns q)$ let us replace the two adjacent
factors in the reordered product \(2.3)
$$
f_{pq}(c_p+z_p\ts,c_q+z_q)\ts
f_{p+1,q}(c_{p+1}+z_{p+1}\ts,c_q+z_q)
$$
by
$$
f_{p,p+\ns1}(c_p\ns+\ns z_p\ts,c_{p+\ns1}\ns+\ns z_{p+\ns1})\,
f_{pq}(c_p+z_p\ts,c_q+z_q)\,
f_{p+1,q}(c_{p+1}+z_{p+1}\ts,c_q+z_q)/2\,.
$$
This does not affect the values of the restriction of \(2.3) to $\D_\om$.
But the restriction to $z_p=z_{p+\ns1}$ of the replacement product
is regular at $z_p=z_q$ by Lemma 2.1.

Due to Proposition 2.4 the divisibility of $F_\om$ by
the element $F_\ps$ on the left is equivalent
to the divisibility by the same element on the right.
Let us now prove the divisibility on the left.
In our new ordering of the factors in \(2.3), we have
$$
F_\om(z_1\lcd z_n)=F_\ps(z_1\lcd z_m)\,\cdot\,
\prod_{(p,q)}^{\longrightarrow}\ f_{pq}(c_p+z_p\ts,c_q+z_q)
\Tag{2.3333}
$$
where $(p,q)\succ(m\!-\!1,\ns m)$. By the above argument, in \(2.3333)
the product over the pairs $(p,\ns q)$ may be
replaced by a certain rational function of $z_1\lcd z_n$ with the restriction
to $\D_\om$ regular at $\or$, without affecting the values of the
restriction. Now the divisibility of $F_\om$ by $F_\ps$
on the left follows from the decomposition \(2.3333)
and the definition of the element $F_\ps\in\CS_m$
\enddemos

\noindent
We can now give a relatively short proof
of the main result of this section, cf.~[N]\>.

\proclaim{Theorem 2.8}
Suppose the numbers $p<q$ stand next to each other in the
same row of the tableau $\Om$. Let $r$ be the number at the bottom
of the column of $\Om$ containing $p$. Then the element $F_\om$
is divisible on the left by the product
$$
\prod_{s\ts=\ts p,...,\ts r}^{\longleftarrow}
\biggl(\
\prod_{t\ts=\ts r+1,...,\ts q}^{\longrightarrow}
f_{st}(\ts c_s,\ns c_t\ts)
\biggr)\,.
\Tag{2.4}
$$
\endproclaim

\demo{Proof}
It suffices to prove Theorem 2.8
only in the particular case $q=n$. For $q<n$ the required statement
would then follow by Proposition 2.7 applied to $m=q$.
Furthermore, it suffices to consider only the case $p=1$.
For $p>1$ the required statement
would then follow by Proposition 2.6 applied to $l=p-1$.

Assume that $p=1$ and $q=n$. Then the skew Young diagram $\om$
consists of two columns only. Further, every row of $\om$ but one consists
of a single box. There is also a row of two boxes, which in the tableau $\Om$
are filled with numbers $1$ and $n$. So there is only one box on every
diagonal of $\om$. The function
$F_\om(z_1\lcd z_n)$ is then regular at $\or$ and takes the value
$$
\gather
\prod_{s\ts=\ts 1,...,\ts n-1}^{\longrightarrow}
\biggl(\
\prod_{t\ts=\ts s+1,...,\ts n}^{\longrightarrow}
f_{st}(\ts c_s,\ns c_t\ts)
\biggr)
\ =
\prod_{s\ts=\ts 1,...,\ts r}^{\longleftarrow}
\biggl(\
\prod_{t\ts=\ts r+1,...,\ts n}^{\longrightarrow}
f_{st}(\ts c_s,\ns c_t\ts)
\biggr)\ \times
\\
\prod_{s\ts=\ts 1,...,\ts r-1}^{\longrightarrow}
\biggl(\
\prod_{t\ts=\ts s+1,...,\ts r}^{\longrightarrow}
f_{st}(\ts c_s,\ns c_t\ts)
\biggr)
\ \cdot
\prod_{s\ts=\ts r+1,...,\ts n}^{\longrightarrow}
\biggl(\
\prod_{t\ts=\ts s+1,...,\ts n}^{\longrightarrow}
f_{st}(\ts c_s,\ns c_t\ts)
\biggr)\,;
\endgather
$$
the last equality was obtained by using the relations \(2.1) and \(2.2).
So the element $F_\om$ is indeed divisible on the left
by the product \(2.4) with $p=1$ and $q=n$
\enddemos


\Section{Leading term near the origin}

\noindent
Take any two non-empty skew Young diagrams $\om$ and $\ompr$.
Let $n$ and $n'$ be the numbers of boxes in $\om$ and $\ompr$
respectively. In this section we will introduce
a certain rational function $F_{\om\om'}(u)$
of one complex variable $u$ with the values in the
symmetric group ring $\CS_{n+n'}$.
Let $z$ and $z'$ be any two complex numbers.
In Section~4,
the leading term of the Laurent expansion of the function $F_{\om\om'}(u)$
at the point $u=z\ns-\ns z'$
will determine an intertwining operator of $\YN\]$-modules
$$
V_\om(z)\ox V_{\om'}(z')\longrightarrow V_\om(z)\oxp V_{\om'}(z')\,.
$$
Here the tilde refers to the opposite comultiplication on
the Hopf algebra $\YN$. In the present section we will compute
the leading term of the Laurent expansion of $F_{\om\om'}(u)$
with $\om=\ompr$ near the origin $u=0$.
This result will be used in Section~4.

Let us regard $\CS_n$ as a subalgebra in $\CS_{n+n'}$
via the standard embedding $S_n\to S_{n+n'}$. In particular, $F_\om\in\CS_n$
may be regarded as an element of $\CS_{n+n'}$. Consider also the embedding
$\CS_{n'}\to\CS_{n+n'}$ defined by
$(p\ts q)\mapsto(p\ns+\ns n,\ns q\ns+\ns n)$ for distinct $p,q=1\lcd n'$.
In this section, the image in $\CS_{n+n'}$
of any element $F\in\CS_{n'}$ under the latter
embedding will be denoted by $\Fv$.

Let $c_1\lcd c_n$ and $\cp_1\lcd\cp_{n'}$ be the
sequences
of contents corresponding to the skew diagrams $\om$ and $\ompr$.
Here we order the contents using the column tableaux.
Define the function $F_{\om\om'}(u)$ as the ordered product
$$
F_\om\ts\cdot\!
\prod_{p=1,...,\ts n}^{\longrightarrow}
\biggl(\
\prod_{q=1,...,\ts n'}^{\longrightarrow}
f_{p,q+n}(\ts c_p\ns+\ns u\ts,\cp_q\ts)
\biggr)\cdot\ts\Fv_{\om'}\,.
\Tag{3.0}
$$
Note that the so defined rational function $F_{\om\om'}(u)$
may have poles only at $u\in\ZZ$.

\proclaim{Proposition 3.1}
The function $F_{\om\om'}(u)$ can be
written as either of the products
$$
\gather
\prod_{p=1,...,\ts n}^{\longleftarrow}
\biggl(\
\prod_{q=1,...,\ts n'}^{\longrightarrow}
f_{p,q+n}(\ts c_p\ns+\ns u\ts,\cp_q\ts)
\biggr)\cdot\ts
F_\om\ts\Fv_{\om'}\,,
\Tag{3.00}
\\
F_\om\ts\Fv_{\om'}\,\,\cdot\!
\prod_{p=1,...,\ts n}^{\longrightarrow}
\biggl(\
\prod_{q=1,...,\ts n'}^{\longleftarrow}
f_{p,q+n}(\ts c_p\ns+\ns u\ts,\cp_q\ts)\biggr)\,.
\Tag{3.000}
\endgather
$$
\endproclaim

\demo{Proof}
Using the definition \(2.3)
along with the relations \(2.1) and \(2.2), we obtain the equality
of the rational functions in $z_1\lcd z_n$ and $u$
$$
\gather
F_\om(z_1\lcd z_n)\ts\cdot\!
\prod_{p=1,...,\ts n}^{\longrightarrow}
\biggl(\
\prod_{q=1,...,\ts n'}^{\longrightarrow}
f_{p,q+n}(\ts c_p\ns+\ns z_p\ns +\ns u\ts,\cp_q\ts)
\biggr)\,=
\\
\prod_{p=1,...,\ts n}^{\longleftarrow}
\biggl(\
\prod_{q=1,...,\ts n'}^{\longrightarrow}
f_{p,q+n}(\ts c_p\ns+\ns z_p\ns +\ns u\ts,\cp_q\ts)\biggr)
\cdot
F_\om(z_1\lcd z_n)\,.
\endgather
$$
By the definition of the element $F_\om\in\CS_n$,
this equality implies that the products \(3.0) and \(3.00)
are also equal to each other. Similarly, we can prove the
equality of the products \(3.0) and \(3.000)
\enddemos

\noindent
Now consider the function $F_{\om\om}(u)$.
For any diagram $\omf$
obtained from $\om\subset\ZZ^2$ by a translation,
we have $F_{{\omf\ns\omf}}(u)=F_{\om\om}(u)$.
Our preliminary result on $F_{\om\om}(u)$ is

\proclaim{Proposition 3.2}
The order of the pole of $F_{\om\om}(u)$ at $u=0$ does not exceed~$d(\om)$.
\endproclaim

\demo{Proof}
We will use the expression for $d(\om)$ given by Proposition 1.3.
Our argument will be similar to that of Proposition 2.2,
but more elaborate, cf.~[\ts N,\ts Theorem 4.1\ts]\>.

The function $F_{\om\om}(u)$ is defined as the ordered product
$$
F_\om\ts\cdot\!
\prod_{p=1,...,\ts n}^{\longrightarrow}
\biggl(\
\prod_{q=1,...,\ts n}^{\longrightarrow}
f_{p,q+n}(\ts c_p\ns+\ns u\ts,c_q\ts)
\biggr)\cdot\Fv_\om\,.
\Tag{3.1}
$$
The factor $f_{p,q+n}(\ts c_p\ns+\ns u\ts,c_q\ts)$ in \(3.1)
has a pole at $u=0$ if and only if the numbers $p$ and $q$ stand on
the same diagonal of the tableau $\Om$. This pole is simple. Again,
we will then call the pair $(p,\ns q)$ singular. Any singular pair
$(p,\ns q)$ belongs to one of the following three types\ts:

\roster
\item"{(i)}" the number $p$ is not at the bottom of its column in $\Om$;
\item"{(ii)}" $p$ is at the bottom of its column in $\Om$ but
$q$ is not leftmost in its row\ts;
\item"{(iii)}" $p$ is at the bottom of its column in $\Om$ and
$q$ is leftmost in its row\ts.
\endroster

\noindent
Observe that the total number of singular pairs of type (iii)
is exactly $d(\om)$. Indeed, let $p$ be at the bottom of the $j\]$-th column.
Then $c_p=j-\lap_j$. Consider the diagonal of the tableau $\Om$
containing $p$. Any number $q$ on this diagonal, except maybe for the
first number, has in $\Om$ a neighbour on the left. As usual, here we are
reading the diagonal from top left to bottom right. Suppose $q$ is the
first on the diagonal. The total number of such pairs $(p,q)$ equals the
number of the columns, which is $\ell(\omp)$. But if $q$ has a neighbour
on the left, say in the $i\]$-th row, then $i-\mup_i\ge c_q$.
If we had $i\ns-\ns\mup_i>c_q$, the number $q$ could not be the first on its
diagonal. So we must have $i\ns-\ns\mup_i=c_q=c_p=j\ns-\ns\lap_j$.
Here we also have $i<j$. Thus the total number
of singular pairs (iii) equals
$$
\ell(\omp)-
\#\ts\{\,(i,\ns j)\ |\ \lap_j\ns-\ns j=\mup_i\ns-\ns i\ts,\ i<j\ts\}
\ts=\ts d(\om)\,.
\nopagebreak
$$
We shall prove that when we estimate
the order of pole of \(3.1) at $u=0$ from above,
the singular pairs (i) and (ii) do not count.
Proposition 3.2 will then follow.

\goodbreak
I.\ Using only the relations \(2.2),
we can rewrite the ordered product \(3.1) as
$$
F_\om\cdot\!\!
\prod_{q=1,...,\ts n}^{\longrightarrow}
\biggl(\
\prod_{p=1,...,\ts n}^{\longrightarrow}
f_{p,q+n}(\ts c_p\ns+\ns u\ts,c_q\ts)
\biggr)\cdot\Fv_\om\,.
\Tag{3.2}
$$
Consider any singular pair $(p,\ns q)$ of type (i). In the ordered product
\(3.2) the factor $f_{p,q+n}(\ts c_p\ns+\ns u\ts,c_q\ts)$ is immediately
followed by factor $f_{p+1,q+n}(\ts c_{p+1}\ns+\ns u\ts,c_q\ts)$.
Note that here $c_p=c_q$.
The numbers $p$ and $p\ns+\ns1$ stand in the same column of the tableau $\Om$.
By Theorem 2.5, the element $F_\om$ is divisible on the right by
$$
f_{p,p+1}(\ts c_p,\ns c_{p+1})=
1\!-\ns(p,p\!+\!1)=
f_{p,p+1}(\ts c_p\ns+\ns u,\ns c_{p+1}\ns+\ns u)\,.
$$
Again using the relations \(2.1) and \(2.2) we obtain, that the product of all
factors in \(3.2) preceding $f_{p,q+n}(\ts c_p\ns+\ns u\ts,c_q\ts)$
is also divisible by  $f_{p,p+1}(\ts c_p\ns+\ns u,\ns c_{p+1}\ns+\ns u)$
on the right. But the product
$$
f_{p,p+1}(\ts c_p\ns+\ns u,\ns c_{p+1}\ns+\ns u)\,
f_{p,q+n}(\ts c_p\ns+\ns u\ts,c_q\ts)\,
f_{p+1,q+n}(\ts c_{p+1}\ns+\ns u\ts,c_q\ts)
$$
has no pole at $u=0$ by Lemma 2.1. Thus
when we estimate the order of pole of \(3.2) at $u=0$ from above,
the singular pairs (i) do not count.

II.\ Now consider the diagonal of $\Om$ containing the number $n$.
Let $q_1\lcd q_d$ be all entries of $\Om$ on that diagonal.
We assume that $q_1<\ldots<q_d$ so that $q_d=n$.
The singular pairs $(n,\ns q_2)\lcd (n,\ns q_d)$ are all of type (ii).
The singular pair $(n,\ns q_1)$ is of type (ii) or (iii). Depending on this,
let the index $k$ range over $1\lcd d$ or over $2\lcd d$ respectively. For
every such $k$ let $p_k$ be the number standing in $\Om$ next to the left of
$q_k$. Of course, for $k>1$ we have $p_k=q_{k-1}\ns+\ns1$. However, to include
the possible case $k=1$ in our argument, we will still use the notation $p_k$.

Let $r_k$ be the number standing in the tableau $\Om$ at the bottom of the
column containing $p_k$. The number at the top of the column containing $q_k$
is then $r_k\ns+\ns1$. By Theorem 2.8, the element $F_\om\in\CS_n$
is divisible on the left by the product
$$
\gather
\prod_{s\ts=\ts p_k,...,\ts r_k}^{\longleftarrow}
\biggl(\
\prod_{t\ts=\ts r_k+1,...,\ts q_k}^{\longrightarrow}
f_{st}(\ts c_s,\ns c_t\ts)
\biggr)\ =
\Tag{3.3}
\\
\prod_{t\ts=\ts r_k+1,...,\ts q_k}^{\longrightarrow}
\biggl(\
\prod_{s\ts=\ts p_k,...,\ts r_k}^{\longleftarrow}
f_{st}(\ts c_s,\ns c_t\ts)
\biggr)\,.
\endgather
$$
For distinct indices $k$ the elements of $\CS_n$ defined by \(3.3)
pairwise commute. Let $F$ be the element of $\CS_n$ obtained
by multiplying the elements \(3.3) over all $k$.
We can write $F_\om=F\ts G$ for some $G\in\CS_n$.
We can also write $\Fv_\om=\Fv\Gv$ in $\CS_{2n}$.

Put $r=r_d$ for short. Using the
relations \(2.2) we can rewrite the~product~\(3.2)~as
$$
\gather
F_\om\cdot\!\!
\prod_{q=1,...,\ts n}^{\longrightarrow}
\biggl(\
\prod_{p=1,...,\ts r}^{\longrightarrow}
f_{p,q+n}(\ts c_p\ns+\ns u\ts,c_q\ts)
\biggr)\ \times
\Tag{3.4}
\\
\prod_{q=1,...,\ts n}^{\longrightarrow}
\biggl(\
\prod_{p=r+1,...,\ts n}^{\longrightarrow}
f_{p,q+n}(\ts c_p\ns+\ns u\ts,c_q\ts)
\biggr)\cdot\Fv_\om\,.
\endgather
$$
Denote by $\Q$ the sequence obtained from the sequence $1\lcd n$
by transposing
its segments $p_k\lcd r_k$ and $r_k\!+\!1\lcd q_k$ for every possible
index $k$. Using the relations \(2.1)\ts,\ts\(2.2) and the decomposition
$\Fv_\om=\Fv\Gv$, we can rewrite~\(3.4)~as
$$
\gather
F_\om\cdot\!\!
\prod_{q=1,...,\ts n}^{\longrightarrow}
\biggl(\
\prod_{p=1,...,\ts r}^{\longrightarrow}
f_{p,q+n}(\ts c_p\ns+\ns u\ts,c_q\ts)
\biggr)\ \times
\Tag{3.5}
\\
\Fv\cdot\,
\prod_{q\ts\in\Q}^{\longrightarrow}\
\biggl(\
\prod_{p=r+1,...,\ts n}^{\longrightarrow}
f_{p,q+n}(\ts c_p\ns+\ns u\ts,c_q\ts)
\biggr)\cdot\ts\Gv\,.
\endgather
$$
As well as in the product \(3.2), here for every singular pair $(p,\ns q)$
of type (i) the factor $f_{p,q+n}(\ts c_p\ns+\ns u\ts,c_q\ts)$ is still
followed by factor $f_{p+1,q+n}(\ts c_{p+1}\ns+\ns u\ts,c_q\ts)$.

Let us now take any singular pair $(n,\ns q_k)$ of type (ii). In the second
line of the ordered product \(3.5), consider the sequence of factors
$$
\prod_{q=q_k,\ts p_k}^{\longrightarrow}\
\biggl(\
\prod_{p=r+1,...,\ts n}^{\longrightarrow}
f_{p,q+n}(\ts c_p\ns+\ns u\ts,c_q\ts)
\biggr)\,.
\Tag{3.6}
$$
Note that $c_p=c_q$ in \(3.6) only if $(p,\ns q)=(n,\ns q_k)$.
The rightmost factor in the product \(3.3) is
$
f_{p_kq_k}(\ts c_{p_k}\ns,\ns c_{q_k})\,.
$
Due to \(2.1) and \(2.2), the product of all factors
in \(3.5) preceding the sequence \(3.6)
is divisible by $f_{p_k+n,q_k+n}(\ts c_{p_k}\ns,\ns c_{q_k})$
on the right. But the cumulate product
$$
\gather
f_{p_k+n,q_k+n}(\ts c_{p_k},c_{q_k})
\ \cdot\!
\prod_{q=q_k,\ts p_k}^{\longrightarrow}\
\biggl(\
\prod_{p=r+1,...,\ts n}^{\longrightarrow}
f_{p,q+n}(\ts c_p\ns+\ns u\ts,c_q\ts)
\biggr)\ =
\\
\prod_{q=p_k,\ts q_k}^{\longrightarrow}\
\biggl(\
\prod_{p=r+1,...,\ts n-1}^{\longrightarrow}
f_{p,q+n}(\ts c_p\ns+\ns u\ts,c_q\ts)
\biggr)\ \times
\\
\phantom{\int}
f_{n,p_k+n}(\ts c_n\ns+\ns u\ts,c_{p_k}\ts)\,
f_{n,q_k+n}(\ts c_n\ns+\ns u\ts,c_{q_k}\ts)\,
f_{p_k+n,q_k+n}(\ts c_{p_k},c_{q_k})
\phantom{\int}
\endgather
$$
has no pole at $u=0$ due to Lemma 2.1.
Thus when we estimate the order of pole of \(3.5) at $u=0$ from above,
any singular pair $(n,\ns q_k)$ of type (ii) does not count.

III. By rewriting the product \(3.4) in its second line as \(3.5),
we excluded
from our count all singular pairs $(p,\ns q)$ of type (ii) with $p=n$.
If the skew Young diagram $\om$ consists of only one column,
there is no singular pairs of type (ii)
and the proof finishes on step I already.
Suppose that the diagram $\om$ has more than one column.
Now let $\ps$ be the skew Young diagram obtained from $\om$ by removing
the last non-empty column. Let $m$ be the total number of boxes in $\ps$.
Due to Proposition~2.7, the element $\Fv_\om$ is divisible by
$\Fv_\ps$ on the left. But we also have
$$
\gather
\prod_{q=1,...,\ts n}^{\longrightarrow}
\biggl(\
\prod_{p=r+1,...,\ts n}^{\longrightarrow}
f_{p,q+n}(\ts c_p\ns+\ns u\ts,c_q\ts)
\biggr)\cdot\Fv_\ps\ =\ \Fv_\ps\ \times
\\
\prod_{q\ts=\ts m,...,1,m+1,...,\ts n}^{\longrightarrow}
\biggl(\
\prod_{p=r+1,...,\ts n}^{\longrightarrow}
f_{p,q+n}(\ts c_p\ns+\ns u\ts,c_q\ts)
\biggr)\,,
\endgather
$$
see Proposition 3.1.
Hence the entire product in the second line of \(3.4)
is divisible on the left by $\Fv_\ps$. Now
arguing similarly to step II, we exclude from our count
all singular pairs $(p,\ns q)$ of type (ii) with $p=m$.
Continuing this argument, we exclude from our count all singular
pairs of type (ii). All singular pairs of type (i)
can be still excluded from our count, as it was done on step I
\enddemos

\noindent
For the skew Young diagram $\om$,  let us introduce the rational
function of $u$
$$
h_\om(u)=\left(\ns\frac{-1}{u}\ns\right)^{\ns\!n}\cdot\prod_{1\le p<q\le n}
\,\biggl(1-\frac1{(u\!-\ns c_p\!+\ns c_q)^2}\biggr)\,.
\Tag{3.625}
$$
\noindent
Recall once again, that here the contents $c_1\lcd c_n$ of $\om$ are ordered
by using the column tableau $\Om$. Now consider the skew Young diagram $\omp$
conjugate to $\om$, and the corresponding rational function $h_\omp\ns(u)$.
In general, we do not have the equality
$h_\omp\ns(u)=(-\ns1)^n\cdot h_\om(-u)$.
But the following proposition is always true.

\proclaim{Proposition 3.3}
The leading term of the Laurent expansion at $u=0$
of the function $h_\omp\ns(u)$
coincides with that of the function $(-\ns1)^n\cdot h_\om(-u)$.
\endproclaim

\demo{Proof}
Consider the sequence of numbers obtained by reading the entries of the column
tableau $\Om$ in the natural way, that is by rows downwards,
from left to right in every row. This sequence can be also obtained from the
sequence $1\lcd n$ by a certain permutation $g\in S_n$.
In general, the permutation $g$ is non-trivial.
The sequence of contents corresponding to the diagram
$\omp$ is \text{$-c_{g(1)}\lcd\!\!-\!c_{g(n)}$.} So
$$
(-\ns1)^n\cdot h_\omp\ns(-u)=
\left(\ns\frac{-1}{u}\ns\right)^{\ns\!n}\cdot\prod_{1\le p<q\le n}
\,\biggl(1-\frac1{(u\!-\ns c_{g(p)}\!\ns+\ns c_{g(q)})^2}\biggr)\,.
$$
Dividing the right hand side of this equality by that of \(3.625),
we get the product~of
$$
\biggl(1-\frac1{(u\ns+\ns c_p\ns-\ns c_q)^2}\biggr)\cdot
\biggl(1-\frac1{(u\ns-\ns c_p\ns+\ns c_q)^2}\biggr)^{\!-1}
\Tag{3.626}
$$
over all pairs $(p,\ns q)$ such that $p<q$ but $g^{-1}(p)>g^{-1}(q)$.
For any such pair, the number $q$ stands above and to the right of $p$
in the tableau $\Om$. In particular, then $c_p\ns-\ns c_q\neq\pm1$.
But then the factor \(3.626) is regular at $u=0$, and
takes the value $1$ at that point
\enddemos

\noindent
The next proposition provides another expression for
the function $h_\om(u)$ with the arbitrary skew Young diagram
$\om=\la\ts/\mu$. Consider the infinite product
$$
\prod_i\,\ts
\frac{u\ns-\ns\la_i\ns+\ns\mu_i}u
\,\ts\cdot\,
\prod_{i<j}\
\frac
{(\ts u\ns+\ns\la_i\ns-\ns\mu_j\ns-\ns i\ns+\ns j\ts)\,
(\ts u\ns+\ns\mu_i\ns-\ns\la_j\ns-\ns i\ns+\ns j\ts)\,}
{(\ts u\ns+\ns\la_i\ns-\ns\la_j\ns-\ns i\ns+\ns j\ts)\,
(\ts u\ns+\ns\mu_i\ns-\ns\mu_j\ns-\ns i\ns+\ns j\ts)\,}
\nopagebreak
\Tag{3.65}
$$
where $i\ts,\ns j=1,\ns2,\ts\ldots\ts$.
Here the first fraction equals $1$
for all but a finite number of indices $i$.
The second fraction equals $1$
for all but a finite number of indices $i$ and $j$.
Hence the product \(3.65) determines a rational function of the variable $u$.

\proclaim{Proposition 3.4}
The rational function \(3.65) is equal to $h_\omp\ns(u)\cdot(1\ns-\ns u)^n$.
\endproclaim

\demo{Proof}
For each $i=1,2,\ts\ldots\ts$ denote by $h_i(u)$ the product of the factors
$$
1-\frac1{(u\ns-\ns c_p\ns+\ns c_q)^2}\,=\,
\frac{\ts u\!-\ns c_p\!+\ns c_q\!-\!1\ts}{u\!-\ns c_p\!+\ns c_q}
\cdot
\frac{\ts u\!-\ns c_p\!+\ns c_q\!+\!1\ts}{u\!-\ns c_p\!+\ns c_q}
\Tag{3.66}
$$
over $p<q$ such that the numbers $p$ and $q$ stand
in the $i\]$-th column of the tableau $\Om$. If this column is empty,
we put $h_i(u)=1$. For these numbers $p$ and $q$ we have
$c_p\!-\ns c_q=q\ns-\ns p$. The total number of boxes in this column is
$\lap_i\ns-\ns\mup_i$. So by direct calcultion
$$
h_i(u)\,=\,
\frac{(u\ns-\ns\lap_i\ns+\ns\mup_i)\,u^{\ts\lap_i-\mup_i-1}}
{(u\ns-\ns 1)^{\ts\lap_i-\mup_i}}\ .
$$

Further, for each $i$ put $a_i=i\ns-\ns\mup_i\ns-\ns1$ and
$b_i=i\ns-\ns\lap_i$. The numbers $a_i$ and $b_i$ are the contents
respectively of the top and the bottom boxes of the
$i\]$-th column of the diagram $\om$, if this column is not empty.
If it is empty, $\lap_i=\mup_i$ and $b_i=a_i\ns+\ns1$.

Now for $i<j$ denote by $h_{ij}(u)$ the product of all those
factors \(3.66), where the numbers $p$ and $q$ stand respectively
in the $i\]$-th and the $j\]$-th columns of $\Om$.
If at least one of these two columns is empty, we put $h_{ij}(u)=1$.
By direct calculation,
$$
h_{ij}(u)
\,=\,
\frac{\ts u\!-\ns a_i\!+\ns b_j\!-\!1\ts}{u\!-\ns a_i\!+\ns a_j}
\cdot
\frac{\ts u\!-\ns b_i\!+\ns a_j\!+\!1\ts}{u\!-\ns b_i\!+\ns b_j}\ .
$$

Note that we also have
$\lap_1\ns-\ns\mup_1\ns+\ns\lap_2\ns-\ns\mup_2+\ldots=n$.
So by the definition \(3.625)
$$
\gather
h_\om(u)\cdot(1\ns-\ns u)^n\,=\,
\left(\ns\frac{u\ns-\ns1}{u}\ns\right)^{\ns\!n}\cdot\,\,
\prod_{i}\ h_i(u)\,\cdot\,
\prod_{i<j}\,h_{ij}(u)\,=\,
\\
\prod_i\,\ts
\frac{u\ns-\ns\lap_i\ns+\ns\mup_i}u
\,\ts\cdot\,
\prod_{i<j}\
\frac
{(\ts u\ns+\ns\mup_i\ns-\ns\lap_j\ns-\ns i\ns+\ns j\ts)\,
(\ts u\ns+\ns\lap_i\ns-\ns\mup_j\ns-\ns i\ns+\ns j\ts)\,}
{(\ts u\ns+\ns\mup_i\ns-\ns\mup_j\ns-\ns i\ns+\ns j\ts)\,
(\ts u\ns+\ns\lap_i\ns-\ns\lap_j\ns-\ns i\ns+\ns j\ts)\,}\ .
\endgather
$$
Replacing here the diagram $\om$ by its conjugate $\omp$, we get
Proposition 3.4
\enddemos

\noindent
The total number of indices $i$ such that $\la_i\neq\mu_i$,
is the number $\ell(\om)$ of non-empty rows in the skew Young diagram $\om$.
For any $i<j$ we also have the inequalities
$$
\la_i\ns-\ns i>\mu_j\ns-\ns j\,,
\quad
\la_i\ns-\ns i>\la_j\ns-\ns j\,,
\quad
\mu_i\ns-\ns i>\mu_j\ns-\ns j\,.
$$
Therefore,
the order of the pole of the rational function \(3.65) at $u=0$ is
$$
\ell(\om)-\#\ts\{\,(i,\ns j)\ |\ \la_j\ns-\ns j=\mu_i\ns-\ns i\ts,\ i<j\ts\}
\ts=\ts d(\om)\ts=\ts d(\omp)\,,
$$
see Proposition 1.2. Now Proposition 3.4 implies
that the order of the pole of the function $h_\om(u)$ at $u=0$ is also
$d(\om)$. Denote by $c(\om)$ the coefficient at $u^{-d(\om)}$
in the Laurent expansion of $h_\om(u)$ at $u=0$.
By Proposition 3.3, we have
$$
c(\omp)\,=\,\ts c(\om)\cdot(-1)^{\ts n+d(\om)}\,.
\Tag{3.777}
$$
Here is the main result of this section.
We will employ this result in Section 4.

\proclaim{Theorem 3.5}
The leading term of the Laurent expansion at $u=0$ of $F_{\om\om}(u)$ is
$$
c(\om)
\cdot
(1,\ns n\ns+\ns1)\,(2,\ns n\ns+\ns2)\,\ldots\,(n,\ns 2n)
\cdot
F_\om\ts\Fv_\om\,
u^{-d(\om)}\,.
\Tag{3.7}
$$
\endproclaim

\demo{Proof}
We will use induction on $n$, the total number of boxes in $\om$. If $n=1$,
then
$$
F_\om=1,\quad
\Fv_\om=1\quad
\text{and}\quad
F_{\om\om}(u)=1\ns-(12)/u\,.
$$
But for $n=1$ we have $h_\om\ns(u)=-1/u$ by the definition \(3.625).
So the statement of Theorem~3.5 is obvious in this case.
{}From now on we will assume that $n>1$.

Suppose that the diagram $\om$ consists of one column only.
Then the diagram $\omp$ consists of only one row.
It follows from Proposition 2.3 and Theorems~2.5\ts,\ts2.7 that then
$$
F_\om=\sum_{g\in S_n}\sgn(g)\cdot g
\ \quad\text{and}\ \quad
F_{\omp}=\sum_{g\in S_n}g
$$
where $\sgn(g)=\pm1$ depending on whether the permutation $g$ is even or odd.
Now denote by $\be$ the involutive automorphism
of the group ring $\CS_{2n}$ defined by $\be(g)=\sgn(g)\cdot g$
for every $g\in S_{2n}$. Then we have
$\be\bigl(F_{\om\om}(u)\bigr)=F_{\omp\ns\omp\ns}(-u)$
by the definition \(3.1). On the other hand, the image
of \(3.7) under $\be$ is
$$
(-1)^n\,c(\om)
\cdot
(1,\ns n\ns+\ns1)\,(2,\ns n\ns+\ns2)\,\ldots\,(n,\ns 2n)
\cdot
F_{\omp\ns}\ts\Fv_{\omp\ns}\,\ts
u^{-d(\om)}\,.
$$
Using the equality \(3.777) we see, that
the statements
of Theorem~3.5 for $\om$ and $\omp$ are equivalent in the case,
when $\om$ consists of one column only. Note that then
$d(\om)=1$, but we did not use the last equality here.
We only used $F_{\omp}=\be(F_\om)$.

{}From now on we will assume that there is more than one column
in the skew Young diagram $\om$.
Let $\de$ be the skew Young diagram consisting
of the first non-empty column of $\om$.
Let $\ph$ the diagram obtained from $\om$
by removing this column, the diagram $\ph$ is also non-empty.
Each of the diagrams
$\de$ and $\ph$ has less than $n$ boxes.
We will assume that the statement of Theorem 3.5 is valid for
each of the diagrams $\de$ and $\ph$ instead~of~$\om$.

I. Let $l$ be the number of boxes in the column $\de$.
As well as in Section 2, denote by $\si_l$ the embedding
$\CS_{n-l}\to\CS_n$ defined by
$\si_l:(p\ts q)\mapsto(p\ns+\ns l,\ns q\ns+\ns l)$
for all distinct $p,q=1\lcd n-l$.
Using \(3.1) with \(2.2), we can write $F_{\om\om}(u)$ as
$$
F_\om\cdot A(u)\ts B(u)\ts C(u)\ts D(u)\cdot\Fv_\om
\Tag{3.8}
$$
where $A(u),B(u),C(u),D(u)$ are ordered products of
the factors $f_{p,q+n}(\ts c_p\ns+\ns u\ts,c_q\ts)$ respectively over
\vskip-27pt
$$
\gather
p,\ns q=1\lcd l\ts;
\\
p=1\lcd l
\quad\text{and}\quad
q=l\ns+\!1\lcd n\ts;
\\
p=l\ns+\!1\lcd n
\quad\text{and}\quad
q=1\lcd l\ts;
\\
p,\ns q=l\ns+\!1\lcd n\ts.
\endgather
$$
The pairs $(p,\ns q)$ in each of these four products are ordered
lexicographically. By Propositions 2.6 and 2.7, the element $F_\om\in\CS_n$
is divisible on the left and on the right by $\si_l(F_\ph)$ and $F_\de$.
The element $F_\de/\ts l\ts!\in\CS_l$ is an idempotent.
\vvgood
In general, there maybe no non-zero idempotents in the subset
$\CC\hskip-.5pt\cdot\hskip-1pt F_\ph\subset\CS_{n-l}$.
But there is always an element $I_\ph\in\CS_{n-l}$
divisible on the left by $F_\ph$, such that
$I_\ph F_\ph=F_\ph$. By Proposition 2.4 applied to $F_\ph$,
this follows from
the semisimplicity of $\CS_{n-l}$. Due to Proposition 2.3,
the element $J_\ph=\al(I_\ph)$ is divisible by $F_\ph$
on the right, and we have $F_\ph\ts J_\ph=F_\ph$.
Using Proposition 3.1, we can now rewrite~\(3.8)~as
$$
F_\om\cdot
F_\de\ts A(u)\ts\Fv_\de\cdot
F_\de\ts B(u)\ts\si_l(I_\ph)^\vee\cdot
\si_l(J_\ph)\ts C(u)\ts D(u)\ts\Fv_\om\cdot(l\ts!)^{-3}\,.
\Tag{3.9}
$$

II. By the inductive assumption for $\de$, the factor
$F_\de\ts A(u)\ts\Fv_\de$ in \(3.9) has a simple pole at $u=0$,
and the corresponding residue coincides with that of
$$
h_\de(u)\cdot (1,\ns n\ns+\ns1)\ldots(l,\ns n\ns+\ns l)
\cdot
F_\de\ts\Fv_\de\,.
\Tag{3.10}
$$
Further, $F_\de\ts B(u)\ts\si_l(F_\ph)^\vee$ is conjugate by the involution
$(l+1,\ns n+l+1)\ldots(n,\ns 2n)$~to
$$
F(u)\,=\,F_\de\ts\cdot\!
\prod_{p=1,...,\ts l}^{\longrightarrow}
\biggl(\
\prod_{q=l+1,...,\ts n}^{\longrightarrow}
f_{pq}(\ts c_p\ns+\ns u\ts,c_q\ts)
\biggr)\cdot\si_l(F_\ph)\,.
\Tag{3.11}
$$
The function $F(u)$ is regular at $u=0$, and $F(0)=F_\om$.
This follows by setting $z_1=\ldots =z_l=u$ in Proposition 2.2.
In particular, the factor $F_\de\ts B(u)\ts\si_l(I_\ph)^\vee$ in \(3.9)
is regular at $u=0$. We will show that the order of the pole at $u=0$
of the factor $\si_l(J_\ph)\ts C(u)\ts D(u)\ts\Fv_\om$ in \(3.9)
does not exceed $d(\om)\!-\!1$. This will allow us to replace
the factor $F_\de\ts A(u)\ts\Fv_\de$ in \(3.9) by the product \(3.10),
without affecting the coefficient at $u^{-d(\om)}$ in the Laurent expansion
at $u=0$. Independently of Proposition 3.2, this will also imply
that the order of the pole at $u=0$ of $F_{\om\om}(u)$
does not exceed $d(\om)$.

It is the argument of Proposition 3.2 that we will use here.
In particular, for any
numbers $p$ and $q$ standing on the same diagonal of $\Om$,
let us keep calling the pair $(p,\ns q)$ singular. Any singular pair belongs
to one of the types (i,ii,iii)
as described in the beginning of the proof
of Proposition 3.2. Note that the number $l$ stands at the bottom
of the first column in the tableau $\Om$. Consider the product
$$
\si_l(J_\ph)\ts C(u)\ts D(u)\ts\Fv_\om\,=\,
\si_l(J_\ph)
\ts\cdot\!
\prod_{p=l+1,...,\ts n}^{\longrightarrow}
\biggl(\
\prod_{q=1,...,\ts n}^{\longrightarrow}
f_{p,q+n}(\ts c_p\ns+\ns u\ts,c_q\ts)
\biggr)\cdot\Fv_\om\,.
$$
Here the factor $f_{p,q+n}(\ts c_p\ns+\ns u\ts,c_q\ts)$
has a pole at $u=0$ if and only if the pair $(p,\ns q)$ is singular.
The factor $\si_l(J_\ph)$ is divisible on the right by $\si_l(F_\ph)$.
By the argument of Proposition 3.2, when we
estimate the order of pole of the above product
at $u=0$ from above,
the singular pairs of types (i) and (ii) do not count.
The total number of the singular pairs $(p,\ns q)$ of type (iii) is $d(\om)$.
For any singular pair $(p,\ns q)$ of type (iii) we have $p>l$,
except for the pair $(p,\ns q)=(l,\ns l)$.
Hence in the above product, the number of singular pairs $(p,\ns q)$
of type (iii) is exactly $d(\om)\!-\!1$.

III. Thus the coefficient at $u^{-d(\om)}$ in the Laurent expansion
at $u=0$ of the function $F_{\om\om}(u)$ coincides with that of the function
$$
\gather
h_\de(u)\,
F_\om\cdot
(1,\ns n\ns+\ns1)\ldots(l,\ns n\ns+\ns l)
\cdot
F_\de\ts\Fv_\de\ \times
\Tag{3.12}
\\
F_\de\ts B(u)\ts\si_l(I_\ph)^\vee\cdot
\si_l(J_\ph)\ts C(u)\ts D(u)\ts\Fv_\om\cdot(l\ts!)^{-3}\ =
\\
h_\de(u)\,F_\om\cdot
(1,\ns n\ns+\ns1)\ldots(l,\ns n\ns+\ns l)
\cdot
F_\de\ts B(u)\ts\si_l(I_\ph)^\vee\cdot
C(u)\ts D(u)\ts\Fv_\om\cdot(l\ts!)^{-1}\,.
\endgather
$$
Without affecting that coefficient, we can also replace in \(3.12)
the constant factor $F_\om$ by the function $F(-u)$; see 
\(3.11). The function $F(-u)$ can be also written~as
$$
\si_l(F_\ph)\ts\cdot\!
\prod_{p=1,...,\ts l}^{\longleftarrow}
\biggl(\
\prod_{q=l+1,...,\ts n}^{\longleftarrow}
f_{pq}(\ts c_p\ts,c_q\ns+\ns u\ts)
\biggr)\cdot F_\de\,,
$$
see Proposition 3.1. By exchanging the indices $p,\ns q$ and
using the commutation relations \(2.2), the last product can be rewritten as
$$
\si_l(F_\ph)\ts\cdot\!
\prod_{p=l+1,...,\ts n}^{\longleftarrow}
\biggl(\
\prod_{q=1,...,\ts l}^{\longleftarrow}
f_{qp}(\ts c_q\ts,c_p\ns+\ns u\ts)
\biggr)\cdot F_\de\,.
\Tag{3.13}
$$

Let us now denote by $E(u)$ the ordered product of
the factors $f_{q+n,p}(\ts c_q\ts,c_p\ns+\ns u\ts)$ over
$$
p=l\ns+\!1\lcd n
\quad\text{and}\quad
q=1\lcd l\,;
$$
the pairs $(p,\ns q)$ being taken in the reversed lexicographical order.
Using the relations \(2.25) repeatedly, we obtain the equality
$$
E(u)\ts C(u)=
\prod_{p=l+1}^n\,
\prod_{q=1}^l
\ \,\biggl(1-\frac1
{(u\!+\ns c_p\!-\ns c_q)^2}
\biggr)\,.
\Tag{3.14}
$$
Denote by $h(u)$ the rational function of $u$ at the right hand side
of this equality.

The product \(3.13) is conjugate by the involution
$(1,\ns n\ns+\ns1)\ldots(l,\ns n\ns+\ns l)\in S_{2n}$ to
$\si_l(F_\ph)\ts E(u)\ts\Fv_\de$. By the above argument,
the coefficient at $u^{-d(\om)}$ in the Laurent expansion
at $u=0$ of the function \(3.12) coincides with that of the function
$$
\gather
h_\de(u)\,F(-u)\cdot
(1,\ns n\ns+\ns1)\ldots(l,\ns n\ns+\ns l)
\cdot
F_\de\ts B(u)\ts\si_l(I_\ph)^\vee\cdot
C(u)\ts D(u)\ts\Fv_\om\cdot(l\ts!)^{-1}\ =
\\
h_\de(u)
\cdot
(1,\ns n\ns+\ns1)\ldots(l,\ns n\ns+\ns l)
\cdot
F_\de\ts B(u)\ts\si_l(I_\ph)^\vee\cdot
\si_l(F_\ph)\ts E(u)\ts\Fv_\de\ts C(u)\ts D(u)\ts\Fv_\om
\cdot(l\ts!)^{-1}
\\
=\ h_\de(u)\ts h(u)
\cdot
(1,\ns n\ns+\ns1)\ldots(l,\ns n\ns+\ns l)
\cdot
F_\de\ts B(u)\ts\si_l(I_\ph)^\vee\cdot
\si_l(F_\ph)\ts D(u)\ts\Fv_\om\,;
\endgather
$$
here we also used divisibility of the product of $\Fv_\om$
and of the product $C(u)\ts D(u)\ts\Fv_\om$ by
$\Fv_\de$ on the left, along with the equality \(3.14).

IV. The element $\Fv_\om$ is also divisible on the left by
$\si_l(F_\ph)^\vee$.
By the inductive assumption for the diagram $\ph$, the leading term in the
Laurent expansion at $u=0$ of
$\si_l(F_\ph)\ts D(u)\ts\si_l(F_\ph)^\vee$ coincides with that of
the function
$$
h_\ph(u)\cdot
(l\ns+\ns1,\ns n\ns+\ns l\ns +\ns\ns1)\ldots(n,\ns 2n)
\cdot
\si_l(F_\ph)\,\si_l(F_\ph)^\vee\,.
$$
Therefore,
by the result of step III,
the coefficient at $u^{-d(\om)}$ in the Laurent expansion
at $u=0$ of the function $F_{\om\om}(u)$ coincides with that of the function
$$
\gather
h_\de(u)\ts h(u)\ts h_\ph(u)\cdot
(1,\ns n\ns+\ns1)\ldots(l,\ns n\ns+\ns l)\ \times
\\
F_\de\ts B(u)\ts\si_l(I_\ph)^\vee\cdot
(l\ns+\ns1,\ns n\ns+\ns l\ns +\ns\ns1)\ldots(n,\ns 2n)\cdot
\si_l(F_\ph)\ts \Fv_\om\ =
\\
h_\de(u)\ts h(u)\ts h_\ph(u)
\cdot
(1,\ns n\ns+\ns1)\ldots(l,\ns n\ns+\ns l)\cdot
(l\ns+\ns1,\ns n\ns+\ns l\ns+\ns\ns1)\ldots(n,\ns 2n)\cdot
F(u)\ts \Fv_\om\,;
\endgather
$$
see the definition \(3.11). But here we have
$h_\de(u)\ts h(u)\ts h_\ph(u)=h_\om(u)$ by \(3.625).
The statement of Theorem 3.5 for the skew Young diagram $\om$ now follows
from the equality $F(0)=F_\om$. This equality
has been already obtained on step II
\enddemos

\noindent
We will complete this section with one remark on the ordered product
$$
\prod_{p=1,...,\ts n}^{\longleftarrow}
\biggl(\
\prod_{q=1,...,\ts n'}^{\longrightarrow}
f_{p,q+n}(\ts c_p\ns+\ns u\ts,\cp_q\ts)
\biggr)
\Tag{3.a}
$$
for two non-empty skew Young diagrams $\om$ and $\ompr$.
This product appeared in \(3.00).
Denote the product \(3.a) by $G_{\om\om'}(u)$, this is a rational function of
$u$ which may have poles only at $u\in\ZZ$. By Proposition 3.1, we have
$F_{\om\om'}(u)=G_{\om\om'}(u)\ts F_\om\ts\Fv_{\om'}$.
Let us now extend the notation $\Fv$ used earlier in this section
for the elements $G\in\CS_{n'}$ as follows.
For any $G\in\CS_{n+n'}$ denote by $G^\vee$
the element $g\ts G g^{-1}$ conjugate to $G$
by the permutation from $S_{n+n'}$
$$
g:\ts
(1\lcd n',n'\!+\!1\lcd n\ns+\ns n')
\,\mapsto\ts
(n\ns+\ns1\lcd n\ns+\ns n',1\lcd n)\,.
$$

\proclaim{Proposition~3.6}
We have the equality
$$
G_{\om\om'}(u)\,G_{\om'\ns\om}^{\ts\vee}(-u)\ts=\ts
\prod_{p=1}^n\,
\prod_{q=1}^{\,n'}
\,\,\biggl(1-\frac1
{(u\!+\ns c_p\!-\ns\cp_q)^2}
\biggr)\,.
$$
\endproclaim

\demo{Proof}
By the definition of the function $G_{\om'\ns\om}(u)$, we have
$$
\align
G_{\om'\ns\om}^{\ts\vee}(-u)
&\ =
\prod_{p=1,...,\ts n'}^{\longleftarrow}
\biggl(\
\prod_{q=1,...,\ts n}^{\longrightarrow}
f_{p+n,q}(\ts\cp_p\ns-\ns u\ts,c_q\ts)
\biggr)
\\
&\ =
\prod_{p=1,...,\ts n}^{\longrightarrow}
\biggl(\
\prod_{q=1,...,\ts n'}^{\longleftarrow}
f_{q+n,p}(\ts\cp_q\ts,c_p\ns+\ns u\ts)
\biggr)\,;
\endalign
$$
here the second equality has been obtained by exchanging the indices
$p\ts,\ns q$ and using the relations \(2.2). Now the required
statement can be derived from the definition \(3.a) of
the function $G_{\om\om'}(u)$
by using the relations \(2.25)
\enddemos


\Section{Irreducibility of Yangian modules}

\noindent
We will begin this section with recalling the definition of the
Yangian of the general linear Lie algebra $\glN$. This is the
unital associative algebra $\YN$ over $\CC$,
with the infinite family of generators
$T_{ij}^{(s)}$ where $s=1,2,\ts\ldots\ts$ and $i\ts,\ns j=1\lcd N$.
The defining relations can be written in terms of the generating
series (0.1) as
$$
(u-v)\cdot[\ts T_{ij}(u)\ts,T_{kl}(v)\ts]\ts=\;
T_{kj}(u)\ts T_{il}(v)-T_{kj}(v)\ts T_{il}(u)\,.
$$
The Yangian $\YN$ is a Hopf algebra.
Coproduct $\De:\YN\to\YN\ox\YN$ is given by
$$
\De\bigl(T_{ij}(u)\bigr)\ts=\ts\sum_{k=1}^N\ T_{ik}(u)\ox T_{kj}(u)\,.
$$
The opposite coproduct $\Det:\YN\to\YN\ox\YN$
is the composition of the coproduct $\De$ and
subsequent transposition of tensor factors in $\YN\ox\YN$.
For any two $\YN\]$-modules $U$ and $V$ we will
denote by $U\ox V$ and $U\oxp V$ the
Yangian modules, where the action of $\YN$ on the tensor
product of vector~spaces is determined via
by the coproduct $\De$ and the opposite coproduct $\Det$ respectively.

Our definition implies, that for any $z\in\CC$
the assignment $T_{ij}(u)\mapsto T_{ij}(u-z)$ defines an automorphism
of the Hopf algebra $\YN$. Here the formal series in $(u-z)^{-1}$ should be
re-expanded in $u^{-1}\]$. Denote this automorphism by $\tau_z$.

We will also need a matrix form of the definition of the Hopf algebra $\YN$.
Let us introduce the following notation.
Consider the tensor product of any unital associative algebras
$\A_1\lox\A_{\ts n}$.
For $p=1\lcd n$ let \text{$\io_p:\A_p\to\A_1\lox\A_{\ts n}$} be the embedding
as the $p\]$-th tensor factor\ts:
$$
\io_p(X)=1^{\ox(p-\ns1)}\ox X\ox1^{\ox(n-p)},
\ \quad
X\in\A_p\,.
$$
For $X\in\A_p$ and $Y\in\A_q$
where $p\neq q$, let us write $(X\ns\ox\ns Y)^{(pq)}=\io_p(X)\ts\io_q(Y)$.
For any $Z\in\A_p\ox\A_q$, determine the element
$Z^{(pq)\!}\in\A_1\ns\lox\ns\A_{\ts n}$ by linearity.
For any vector spaces $U_1\lcd U_n$ we identify the algebras
\text{$\End(U_1\lox U_n)$} and $\End(U_1)\lox\End(U_n)$.

Now take the vector space $\CC^N=\U$. For $i\ts,j=1\lcd N$
let $E_{ij}\in\End(\U)$ be the matrix units.
Let $P:x\ox y\mapsto y\ox x$ be the flip map
in $\U\ox\U$. Put $R(u)=u-\ns P$.
Identifying $\End(U\ox U)=\End(U)\ox\End(U)$, we can write
$$
R(u)\ts=\,u\,-\sum_{i,j=1}^N E_{ij}\ox E_{ji}\,.
\Tag{R}
$$
Combine all series (0.1) into a series $T(u)$ with coefficients
in $\YN\ox\End(\U)$:
$$
T(u)=\sum_{i,j=1}^N T_{ij}(u)\ox E_{ji}\,.
\Tag{T}
$$
Then the defining relations of the algebra $\YN$ are equivalent
to the relation
$$
T^{(01)}(u)\,T^{(02)}(v)\,R^{(12)}(u-v)
\,=\,
R^{(12)}(u-v)\,T^{(02)}(v)\>T^{(01)}(u)
\Tag{RTT}
$$
for the series with coefficients in $\YN\ox\End(\U)\ox\End(\U)$,
where we label the tensor factors by the indices $0\ts,\ns1\ts,2$.

Now label the tensor factors in $\YN\ox\YN\ox\End(\U)$
by $1,\ns2,\ns3$ as usual.~Put
$$
T^{[k]}(u)=T^{(k3)}(u)\,;
\ \quad
k=1\ts,\ns2.
\Tag{ai}
$$
Then we have the equality
$$
(\De\ox\id)\bigl(\ts T(u)\bigr)=T^{[2]}(u)\,T^{[1]}(u)\,.
\Tag{DlT}
$$
We refer the reader to the survey [MNO] for more details on the definition of
the Yangian $\YN$, and for description of its basic properties.
Note that the series denoted by $T(u)$ in [MNO]
differs from our series \(T) by transposition $E_{ij}\mapsto E_{ji}$
in the second tensor factor.

Let $\om=\la\ts/\mu$ be a non-empty skew Young diagram.
Let $n$ be the number of boxes in $\om$.
Take the vector space $\Ux n\]$.
The symmetric group $S_n$ acts naturally in this vector space,
permuting the tensor factors. Denote by $\pi_n$
the corresponding homomorphism $\CS_n\to\End(\Ux n)$.
In Section~2 we defined a certain element $F_\om\in\CS_n$. Let the
subspace $V_\om\subset\Ux n$ be the image of the operator~$\pi_n(F_\om)$.

\proclaim\nofrills{Proposition~4.1}{\bf\,:}
$V_\om\neq\{0\}$ if and only if
the length of every column of $\om$ does~not exceed $N$.
\endproclaim

\noindent
The results of the present article do not depend on this proposition.
Hence we omit the proof, and refer to \text{[C1\ts,\ts NT2\ts]} instead.

Now take two non-empty skew Young diagrams $\om$ and $\ompr$ with $n$
and $n'$ boxes respectively.
Let $c_1\lcd c_n$ and $\cp_1\lcd\cp_{n'}$ be the corresponding
sequences of contents,
ordered according to the column tableaux.
Consider the rational function $G_{\om\om'}(u)$ with the
values in $\CS_{n+n'}$, defined by \(3.a).
For any value of $u\notin\ZZ$
the operator $\pi_{n+n'}(G_{\om\om'}(u))$ is well defined,
and by Proposition~3.1 preserves the subspace
$$
V_\om\ox V_{\om'}\subset\Ux n\!\ox\Ux{n'}\!\!=\Ux{(n+n')}.
$$
Let $R_{\ts\om\om'}(u)$ be the value of the restriction of
$\pi_{n+n'}(G_{\om\om'}(u))$ to $V_\om\ox V_{\om'}$;
$R_{\ts\om\om'}(u)$ is a rational function of $u$ valued in
$\End(V_\om)\ox\End(V_{\om'})$. Suppose that each of the vector
spaces $V_\om$ and $V_{\om'}$ is non-zero, see Proposition~4.1.
Then by Proposition~3.6
$$
R_{\ts\om\om'}(u)\,R_{\om'\ns\om}^{(21)}(-u)\ts=\ts
\prod_{p=1}^n\,\prod_{q=1}^{\,n'}\,
\Bigl(\>1-\frac1{(u+c_p-\cp_q)^2}\>\Bigr)\,.
\Tag{RR}
$$
Take one more non-empty skew Young diagram $\om^{\ts\prime\prime}$.
Let $n''$ be its number~of~boxes.

\proclaim{Proposition~4.2}
We have the equality of rational functions of $u\ts,\ns v$ with
values in the tensor product\/
$\End(V_{\om})\ox\End(V_{\om'})\ox\End(V_{\om''})$
$$
R_{\om\om'}^{(12)}(u-v)\,R_{\om\om''}^{(13)}(u)\,R_{\om'\ns\om''}^{(23)}(v)
\ts=\ts
R_{\om'\ns\om''}^{(23)}(v)\,R_{\om\om''}^{(13)}(u)\,R_{\om\om'}^{(12)}(u-v)\,.
\Tag{RRR}
$$
\endproclaim

\demo{Proof}
The relation \(RRR) for the functions
$R_{\om\om'}(u)$,$R_{\om\om''}(u)$,$R_{\om'\ns\om''}(u)$
follows from the respective relation for the functions
$G_{\om\om'}(u)$,$G_{\om\om''}(u)$,$G_{\om'\ns\om''}(u)\ns$.
The latter relation is an equality of rational functions of $u\ts,\ns v$
with the values in the algebra ${\CS_{n+n'+n^{\prime\prime}}}$.
That equality follows from \(2.1) and \(2.2), we will not write it here.
By applying the homomorphism $\pi_{n+n'+n^{\prime\prime}}$ to both
sides of that equality and taking the restriction to
$V_{\om}\ox V_{\om'}\ox V_{\om''}\subset\Ux{(n+n'+n'')}$,
we obtain \(RRR)
\enddemos

\noindent
Denote by $\ep$ the Young diagram with one box $(1,\ns1)$.
Then $V_\ep=U$. Note that $R_{\ts\ep\ep}(u)=R(u)/u$.
We will now define a family of modules over the algebra $\YN$,
we call these modules elementary. Take any $z\in\CC$.
Consider the rational function $R_{\om\ep}(z\ns-\ns u)$ of the variable $u$.
This function takes values in $\End(V_\om)\ox\End(\U)$.
This function is regular at $u=\infty$. Moreover, at $u=\infty$
it takes the value $1$. Here we keep to the assumption $V_\om\neq\{0\}$.
Consider the Laurent
expansion at $u=\infty$
$$
R_{\ts\om\ep}(z-u)\ts=\ts
1\,+\sum_{i,j=1}^N u^{-s}\ts\,E_{ij}^{(s)}(\om,\ns z)\ox E_{ji}\ts\,.
$$
Here each tensor factor $E_{ij}^{(s)}(\om,\ns z)$ is a certain element of
the algebra $\End(V_\om)$.

\proclaim{Proposition~4.3}
The assignment for all\/ $s=1,2,\ts\ldots\ts$ and\/ $i\ts,\ns j=1\lcd N$
$$
T_{ij}^{(s)}\mapsto\ts E_{ij}^{(s)}(\om,\ns z)
\Tag{Tt}
$$
defines an $\YN\]$-module structure on the vector space $V_\om$.
\endproclaim

\demo{Proof}
We have to verify that the operators $E_{ij}^{(s)}(\om,\ns z)$
satisfy the defining relations of $\YN$.
Let us rewrite the definition \(Tt) in a matrix form. It will be
$$
T(u)\mapsto\ts R_{\ts\om\ep}(z\ns -\ns u)\,.
\Tag{TR}
$$
In view of the defining relation \(RTT), we have to verify the equality
$$
R_{\om\ep}^{(01)}(z\ns-\ns u)\,
R_{\om\ep}^{(02)}(z\ns-\ns v)\,
R^{(12)}(u\ns-\ns v)\ts=\ts
R^{(12)}(u\ns-\ns v)\,
R_{\om\ep}^{(02)}(z\ns-\ns v)\,
R_{\om\ep}^{(01)}(z\ns-\ns u)
$$
of rational functions in $u\ts,\ns v$ with values in
$\End(V_{\om})\ox\End(U)\ox\End(U)$, where the tensor factors
are labelled by the indices $0\ts,\ns1\ts,2$. Changing these labels
respectively to $1\ts,\ns2\ts,3$ and taking into account that
$R(u)=u\ts R_{\ts\ep\ep}(u)$, we come to verifying
$$
R_{\om\ep}^{(12)}(z\ns-\ns u)\,
R_{\om\ep}^{(13)}(z\ns-\ns v)\,
R_{\ep\ep}^{(23)}(u\ns-\ns v)\ts=\ts
R_{\ep\ep}^{(23)}(u\ns-\ns v)\,
R_{\om\ep}^{(13)}(z\ns-\ns v)\,
R_{\om\ep}^{(12)}(z\ns-\ns u)\,.
$$
But the latter equality follows from Proposition~4.2 when
$\ompr=\om^{\ts\prime\prime}=\ep$
\enddemos

\noindent
The $\YN\]$-module defined by Proposition~4.3 will be denoted by $V_\om(z)$,
and called an {\it elementary\/} module. Note that for any $z\in\CC$ the
$\YN\]$-module $V_\om(z)$ is obtained from the module $V_\om(0)$ by
pulling back through the automorphism $\tau_z$ of the algebra $\YN$.

Let us now give another description of the elementary
$\YN\]$-module $V_\om(z)$. First consider the module
$V_\ep(z)$. The vector space of this module is $U=\CC^N$, and the action
of the algebra $\YN$ in this space is determined by the assignment
$T_{ij}^{(s)}\mapsto\ts z^s E_{ij}$, see \(R).

\proclaim{Proposition~4.4}
The action of the algebra $\YN$ in the tensor product of the
modules $V_\ep(c_1\ns+\ns z)\lox V_\ep(c_n\ns+\ns z)$ preserves the subspace
$V_\om\subset\Ux n$. The restriction of this action to $V_\om$
gives exactly the $\YN\]$-module $V_\om(z)$.
\endproclaim

\demo{Proof}
Let us write the definition of the action of the algebra
$\YN$ in $V_\ep(z)$
in a matrix form. It will
be $T(u)\mapsto R_{\ts\ep\ep}(z\ns-\ns u)$. Then due to \(DlT) the action of
$\YN$ in the tensor product $V_\ep(c_1\ns+\ns z)\lox V_\ep(c_n\ns+\ns z)$
is determined by the assignment
$$
T(u)\ts\mapsto\prod_{p=1,...,\ts n}^{\longleftarrow}
R_{\ts\ep\ep}^{(p,n+1)}(c_p+z-u)\,.
$$
The right hand side of this assignment can be rewritten as
$$
\prod_{p=1,...,\ts n}^{\longleftarrow}
\pi_{n+1}\bigl(f_{p,n+1}(c_p+z-u\ts,\ns0)\bigr)
\,=\,\ts
\pi_{n+1}\bigl(\ts G_{\om\ep}(z\ns-\ns u)\bigr)\,.
$$
We have already shown that $\pi_{n+1}\bigl(G_{\om\ep}(z-u)\bigr)$ preserves
the subspace $V_\om\ox V_\ep$ in $\Ux n\ox U$. This means that the action of
$\YN$ in $V_\ep(c_1\ns+\ns z)\lox V_\ep(c_n\ns+\ns z)$ preserves
the subspace $V_\om\subset\Ux n$.
Restriction of this action to $V_\om$ is
determined exactly by the assignment
$$
T(u)\ts\mapsto\ts R_{\ts\om\ep}(z\ns-\ns u)
$$
because of the definition of $R_{\ts\om\ep}(z\ns-\ns u)$.
But the last assignment coincides with the matrix
definition of the $\YN\]$-module $V_\om(z)$
\enddemos

\noindent
For any $z\in\CC$ let us expand the rational function $R_{\om\om'}(u)$
into the Laurent series at $u=z$. Let
$(u\ts-z)^{-a_{\om\om'}(z)}\,I_{\om\om'}(z)$
be the leading term of this expansion.
This defines an integer $a_{\om\om'}(z)$ and a non-zero linear
operator $I_{\om\om'}(z)$ in $V_\om\ox V_{\om'}$.
Again, here we keep to the assumption that $V_\om\,,V_{\om'}\neq\{0\}$.
The rational function $R_{\om\om'}(u)$ is regular at $u=z$
for any $z\notin\ZZ$. Then by \(RR), this function also
does not vanish at $u=z$ for any $z\notin\ZZ$. Thus for any $z\notin\ZZ$
we have $a_{\om\om'}(z)=0$ and $I_{\om\om'}(z)=R_{\ts\om\om'}(z)$.
Furthermore, the equality \(RR) implies

\proclaim{Proposition~4.5}
For any $z\in\CC$ the product
$I_{\om\om'}(z)\,I_{\om'\ns\om}^{(21)}(-z)\in\End(V_\om\ox V_{\om'})$ is
proportional to the identity map. In particular, the maps $I_{\om\om'}(z)$
and $I_{\om'\ns\om}(-z)$ are invertible simultaneously.
\endproclaim

\noindent
The basic property of
the linear operator $I_{\om\om'}(z)$ is given by the next proposition.

\proclaim{Proposition~4.6}
For any $z,z'\in\CC$ the element
$I_{\om\om'}(z-z')\in\End(V_\om\ox V_{\om'})$ is an intertwiner
of the Yangian modules
$V_\om(z)\ox V_{\om'}(z')\ts\longrightarrow\ts V_\om(z)\oxp V_{\om'}(z')$.
\endproclaim

\demo{Proof}
The formula \(RRR)
implies the following equality of rational functions in $u,v$:
$$
R_{\om\om'}^{(12)}(z\ns-\ns v)\,
R_{\om\ep}^{(13)}(z\ns-\ns u)\,
R_{\om'\ns\ep}^{(23)}(v\ns-\ns u)
\ts=\ts
R_{\om'\ns\ep}^{(23)}(v\ns-\ns u)\,
R_{\om\ep}^{(13)}(z\ns-\ns u)\,
R_{\om\om'}^{(12)}(z\ns-\ns v)\,.
$$
Taking the leading terms of the Laurent expansions of both sides at $v=z'$
we get
$$
I_{\om\om'}^{\>(12)}(z\ns-\ns z')\,
R_{\om\ep}^{(13)}(z\ns-\ns u)\,
R_{\om'\ns\ep}^{(23)}(z'\ns-\ns u)
\ts=\ts
R_{\om'\ns\ep}^{(23)}(z'\ns-\ns u)\,
R_{\om\ep}^{(13)}(z\ns-\ns u)\,
I_{\om\om'}^{\>(12)}(z\ns-\ns z')\,.
$$
The latter equality implies Proposition~4.6, see \(DlT) and \(TR)
\enddemos

\noindent
Recall that the positive integer
$d(\om)$ is the Durfee rank of the non-empty skew Young diagram
$\om$, see Section~1. The non-zero rational number
$c(\om)$ is the value of the function $u^{\ts d(\om)}\ts h_\om(u)$
at $u=0$, see \(3.625). Denote
by $P_\om$ be the flip map in $V_\om\ox V_\om$, so that
$P_\om(x\ox y)=y\ox x$ for any $x,y\in V_\om$.
The main result of the previous section, Theorem 3.5,
can be now restated as follows.

\proclaim{Proposition~4.7}
We have $a_{\om\om}(0)=d(\om)$ and $I_{\om\om}(0)=c(\om)\ts P_\om$.
\endproclaim

\demo{Proof}
The rational function $F_{\om\om}(u)$ with the values in $\CS_{2n}$
can be written as $G_{\om\om}(u)\ts F_\om\ts\Fv_\om$. Hence the required
statement follows from Theorem~3.5, by the definition of the operator
$I_{\om\om}(0)$ in $V_\om\ox V_\om$
\enddemos

\noindent
Now take $k\ge1$ non-empty skew Young diagrams $\om_1\lcd\om_k$.
We will assume that each of the
vector spaces $V_{\om_1}\lcd V_{\om_k}$
is not zero, see Proposition~4.1. Take also $k$ arbitrary complex numbers
$z_1\lcd z_k$. Consider the elementary $\YN\]$-modules
$V_{\om_1}(z_1)\lcd V_{\om_k}(z_k)$.
The next theorem is the main result of the~present~paper.

\proclaim{Theorem~4.8}
The $\YN\]$-module
$V_{\om_1}(z_1)\lox V_{\om_k}(z_k)$ is irreducible if and only if all
the intertwiners $I_{\om_i\om_j}(z_i\ns-\!z_j)$ with $1\le i<j\le k$ 
are invertible.
\endproclaim

\noindent
There is an immediate corollary, we state
it as a theorem because of its importance.

\proclaim{Theorem~4.9}
The $\YN\]$-module
$V_{\om_1}(z_1)\lox V_{\om_k}(z_k)$ is irreducible if and only if all
the modules $V_{\om_i}(z_i)\ox V_{\om_j}(z_j)$ with $1\le i<j\le k$
are irreducible.
\endproclaim

\noindent
We shall prove Theorem 4.8 in the remainder of this section.
However, before giving the proof, let us make a few remarks.
Applying Theorem~4.8 when $k=1,2$ we obtain, in particular,
the following properties of the elementary $\YN\]$-modules:

\smallskip
\bitem{(a)}
the $\YN\]$-module $V_\om(z)$ is irreducible;
\bitem{(b)}
the $\YN\]$-module $V_\om(z)\ox V_{\om'}(z')$ is irreducible provided
$z-z'\notin\ZZ$;
\bitem{(c)}
if $V_\om(z)\ox V_{\om'}(z')$ is irreducible, then it is equivalent to
$V_\om(z)\oxp V_{\om'}(z')$.

\smallskip
\noindent
The properties (a,b,c) are well known, see for instance
[\ts C2\ts,CP,\ts NT1\ts]\ts.
But we would like to emphasize here,
that Theorem~4.8 gives each of them a new independent proof.
Using  Theorem~4.8 together with Proposition 4.7, we obtain
for any $k=2,3,\ts\ldots$
a new remarkable property of the elementary $\YN\]$-modules:

\smallskip
\bitem{(d)}
the $k\]$-th tensor power of the $\YN\]$-module $V_\om(z)$
is irreducible.

\smallskip
\noindent
Conversely, the property (d) with $k=2$ 
implies that the operator $I_{\om\om}(0)$ in the vector space
$V_\om\ox V_\om$ is proportional to the flip map $P_\om$;
see Proposition 4.6.

The {\it only if\/} part of Theorem 4.8 is an easy and well known statement.
For the sake of completeness, with will give the proof of this statement
at the end of this section. Our proof of the {\it if\/} part is based
on the following general observation. Let $\A$ be any
unital associative algebra over $\CC$,
and let $W$ be any finite-dimensional $\A\]$-module.
Denote by $\rho$ the corresponding homomorphism $\A\to\End(W)$. Let
$P_{\ts W}\in\End(W\ox W)$ be the flip map.

\proclaim{Lemma~4.10}
If $P_{\ts W}\in\rho\ts(\A)\ox\End(W)$,
then the $\A\]$-module $W$ is irreducible.
\endproclaim

\demo{Proof}
Let $1\in\End(W)$ be the unit element, we have $1\in\rho\ts(\A)$
by our assumption. Furthermore, we have
$(1\ox X)\cdot P_{\ts W}\in\rho\ts(\A)\ox\End(W)$ for any $X\in\End(W)$.
Now let $\tr\!:\,\End(W)\to\CC$ be the trace map.
Then the subalgebra $\rho\ts(\A)\subset\End(W)$ contains the element
$$
(\ts\id\ox\tr)\bigl((1\ox X)\cdot P_{\ts W}\ts\bigr)=X\,.
$$
So we actually have $\rho\ts(\A)=\End(W)$
\enddemos

\noindent
Now take the algebra $\A=\YN$ and its module
$W=V_{\om_1}(z_1)\lox V_{\om_k}(z_k)$.
Under the assumption that all
the intertwiners $I_{\om_i\om_j}(z_i\ns-\!z_j)$ with $1\le i<j\le k$
are invertible, we will show that in the above notation,
$$
P_{\ts W}\in\rho\,(\YN)\ox\End(W)\,.
\Tag{PW}
$$
The {\it if\/} part of Theorem~4.8 will then follow by Lemma~4.10.
We will split the proof of \(PW) into four propositions below.

Let $\om$ be a non-empty skew Young diagram with $n$ boxes.
We keep assuming that the subspace $V_\om\subset\Ux n$ is
non-zero. Let
$\End_\om(\Ux n)\subset\End(\Ux n)$ be the stabilizer of $V_\om$. Let
$$
\chi:\,\End_\om(\Ux n)\to\End(V_\om)
$$
be the canonical homomorphism\ts:
$\chi(X)=X|_{V_\om}$.
Introduce the series $\Tti_\om(u)$ with coefficients in
the tensor product \text{$\YN\ox\End(\Ux n)=\YN\ox(\End(\U))^{\ox n}$}
$$
\Tti_\om(u)=T^{(01)}(u+c_1)\,\,\ldots\,\>T^{(0n)}(u+c_n)\,.
$$
Here we labelled the factors in $\YN\ox(\End(\U))^{\ox n}$
by the indices $0\ts,\ns1\lcd n$.

\proclaim{Proposition~4.11}
Coefficients of the series $\Tti_\om(u)$ are in\/ $\YN\ox\End_\om(\Ux n)$.
\endproclaim

\demo{Proof}
Consider the rational function $F_\om(z_1\lcd z_n)$, defined
in Section~2 as the product (2.5). The relation \(RTT) implies that
$$
\gather
T^{(01)}(u\ns+\ns c_1\ns+\ns z_1)\,\,\ldots\,\>
T^{(0n)}(u\ns+\ns c_n\ns+\ns z_n)\cdot
\pi_n\bigl(F_\om(z_1\lcd z_n)\bigr)\ts=
\\
\pi_n\bigl(F_\om(z_1\lcd z_n)\bigr)\cdot
T^{(0n)}(u\ns+\ns c_n\ns+\ns z_n)\,\ldots\,\ts
T^{(01)}(u\ns+\ns c_1\ns+\ns z_1)\,.
\endgather
$$
Hence, by Proposition~2.2
$$
\Tti_\om(u)\,\pi_n(F_\om)\ts=\ts
\pi_n(F_\om)\,T^{(0n)}(u\ns+\ns c_n)\,\ldots\,T^{(01)}(u\ns+\ns c_1)
\quad\square
$$
\enddemo

\noindent
Due to the Proposition~4.11, we can apply the homomorphism $\id\ox\chi$ to
the series $\Tti_\om(u)$, and obtain in this way a series with the
coefficients
in $\YN\ox\End(V_\om)$. We will denote the latter series by $T_\om(u)$.

Henceforth, we will use the following generalization of notation \(ai).
Given a vector space $V$ and an element $X\in\YN\ox\End(V)$, we will write
$$
X^{[k]}=X^{(k3)}\in\YN\ox\YN\ox\End(V)\,,\ \quad k=1\ts,\ns2\,.
$$

\proclaim\nofrills{Proposition~4.12}{\bf\,:}
$\De\bigl(T_\om(u)\bigr)=T^{\ts[2]}_\om(u)\,T^{\ts[1]}_\om(u)$.
\endproclaim

\demo{Proof}
The formula \(DlT) implies that
$\De\bigl(\Tti_\om(u)\bigr)=\Tti^{\ts[2]}_\om(u)\,\Tti^{\ts[1]}_\om(u)$.
By applying the homomorphism $\id\ox\id\ox\chi$ to the last equality,
the proposition follows
\enddemos

\noindent
Let $\ga:\YN\to\End(V_\om)$ be the defining homomorphism
of the module $V_\om(z)$. Take one more skew Young diagram $\ompr$,
we will again assume that $V_{\om'}\neq\{0\}$.

\proclaim\nofrills{Proposition~4.13}{\bf\,:}
$(\ts\ga\ox\id)\bigl(T_{\om'}(u)\bigr)=R_{\ts\om\om'}(z-u)$.
\endproclaim

\demo{Proof}
Let $n'$ be the numbers of boxes in the diagram $\ompr$.
Consider the canonical homomorphism
$$
\chip:\,\End_{\om'}(\Ux{n'})\to\End(V_{\ompr})\,.
$$
It is straighforward to verify that each of
$(\ga\ox\id)\bigl(\ts T_{\om'}(u)\bigr)$ and $R_{\ts\om\om'}(z-u)$
equals
$$
\bigl((\,\chi\ox\chip)\crc\pi_{n+n'}\bigr)
\bigl(G_{\om\om'}(z-u)\bigr)\quad\square
$$
\enddemo

\noindent
Now consider the formal power series in $u^{-1}$
$$
T_W(u)\ts=\ts
T_{\om_1}^{(01)}(u\ns+\ns z_1)\,\,\ldots\,\ts
T_{\om_k}^{(0k)}(u\ns+\ns z_k)
\in\YN\ox\End(W)\,,
$$
here we identify $\YN\ox\End(W)$ with
$\YN\ox\End(V_{\om_1})\lox\End(V_{\om_k})$,
and label the tensor factors by $0,1\lcd k$.
The next proposition completes the proof of
the {\it if\/} part of Theorem 4.8.

\proclaim{Proposition 4.14}
The image $(\ts\rho\:\ox\id)\bigl(T_W(u)\bigr)$
is a rational function in $u$ with values in
$$
\rho\ts(\YN)\ox\End(W)\subset\End(W)\ox\End(W)\,.
$$
The leading coefficient
in its Laurent expansion at $u=0$ is proportional to $P_{\ts W}$.
\endproclaim

\demo{Proof}
It follows from Proposition~4.12 that
$
\De\bigl(T_W(u)\bigr)\ts=\ts T^{\ts[2]}_W(u)\,T^{\ts[1]}_W(u)\,.
$
Using this relation and Proposition~4.13, we get
$$
(\ts\rho\:\ox\id)\bigl(T_W(u)\bigr)\,=
\prod_{j=1,...,\ts k}^{\longrightarrow}
\biggl(\ \prod_{i=1,...,\ts k}^{\longleftarrow}
R_{\ts\om_i\om_j}^{\ts(i,j+k)}(z_i\ns-\ns z_j\ns-\ns u)\biggr)\,.
\Tag{4.99}
$$
Here we regard $\End(W)\ox\End(W)$ as the $2k\]$-fold tensor product
$$
\End(V_{\om_1})\lox\End(V_{\om_k})\ox\End(V_{\om_1})\lox\End(V_{\om_k})\,.
\Tag{4.9}
$$
The right hand side of the equality \(4.99) is a rational function in $u$
by definition. Denote
$$
a\,\ts=\sum_{i,j=1}^k\,a_{\ts\om_i\om_j}(z_i\ns-\ns z_j)\,.
$$
Then the product of the leading terms of the Laurent expansions
at $u=0$ of all factors in the product \(4.99) equals $(-u)^{-a}\]$
multiplied by
$$
\prod_{j=1,...,\ts k}^{\longrightarrow}
\biggl(\ \prod_{i=1,...,\ts k}^{\longleftarrow}
I_{\ts\om_i\om_j}^{\ts(i,j+k)}(z_i\ns-\ns z_j)\biggr)\,.
\Tag{4.999}
$$
By Proposition~4.7, any factor $I_{\ts\om_i\om_j}^{\ts(i,i+k)}(0)$
corresponding to $i=j$, equals~$P_{\om_i}^{\>(i,i+k)}$ up to a non-zero
multiplier. Using this result, we can rewrite \(4.999) as
$$
\prod_{i=1}^k\,P_{\om_i}^{\>(i,i+k)}\,\cdot\!
\prod_{1\le i<j\le k}
\Bigl(
I_{\om_j\om_i}^{\>(j+k,i)}(z_j\ns-\ns z_i)\,
I_{\om_i\om_j}^{\>(i,j+k)}(z_i\ns-\ns z_j)
\Bigr)
\Tag{PII}
$$
up to a nonzero multiplier.
Here all factors $P_{\om_i}^{\>(i,i+k)}$ pairwise commute.
By Proposition~4.5,
each of the factors in \(PII) corresponding to $1\le i<j\le k$, is
proportional to the identity map in \(4.9).
In particular, these factors also pairwise commute.
Moreover, each of these factors is non-zero, because
the operators $I_{\om_i\om_j}(z_i\ns-\ns z_j)$
with $1\le i<j\le k$ are invertible by our assumption.
So the entire product \(4.999) is non-zero
and proportional to
$$
\prod_{i=1}^k\,P_{\om_i}^{\>(i,i+k)}=P_{\ts W}\,.
$$
Thus the leading term in the Laurent expansion
of the function $(\ts\rho\:\ox\id)\bigl(T_W(u)\bigr)$
near the origin $u=0$
equals $(-u)^{-a}\>P_{\ts W}$, up to a non-zero factor from $\CC$
\enddemos

\noindent
We will complete this section with the proof of the {\it only if\/}
part of Theorem~4.8.
Consider the $\YN\]$-module $W=V_{\om_1}(z_1)\lox V_{\om_k}(z_k)$. Suppose
there is a non-invertible
intertwiner $I_{\om_i\om_j}(z_i\ns-\ns z_j)$ with
$1\le i<j\le k$. Assume that the pair $(i\ts,\ns j)$
is lexicographically minimal with this non-invertibility property.
The~operator
$$
I_{\om_i\om_j}^{\>(ij)}(z_i\ns-\ns z_j)\,
\ldots\,
I_{\om_i\om_{i+1}}^{\>(i,i+1)}(z_i-z_{i+1})
$$
is non-zero and non-invertible. By Proposition~4.6,
this is an intertwiner of Yangian modules
\vskip-25pt
$$
\gather
W\ \longrightarrow\
V_{\om_1}(z_1)\lox V_{\om_{i-1}}(z_{i-1})\ \ts\ox
\\
\bigl(\,V_{\om_i}(z_i)\oxp
\bigl(\,V_{\om_{i+1}}(z_{i+1})\lox V_{\om_j}(z_j)\bigr)\bigr)\,\ox\,
V_{\om_{j+1}}(z_{j+1})\lox V_{\om_k}(z_k)\,.
\endgather
$$
So the $\YN\]$-module $W$ is reducible.
This proves the {\it only if\/} part of Theorem~4.8.


\section{Acknowledgements}

\noindent
It is our pleasure to thank N.\,Kitanine for a stimulating discussion
and for bringing the results of [\ts KMT,\ts MT\ts] to our attention.
This work was done when the authors stayed at the
Mathematisches Forschungsinstitut Oberwolfach, with the support
from the Volkswagen-Stiftung. While finishing this article, the second
author stayed at the Max-Planck-Institut f\"ur Mathematik in Bonn.
The first author has been supported by the
EPSRC, and by the EC under \text{the TMR grant FMRX-CT97-0100.}
The second author has been supported by the RFFI grant 99-01-00101 and by
the INTAS grant 99-01705.


\section{References}
\refskips
\widest{[MT]}

\smallskip
\refitem{[A]}
{G.\,Andrews},
\emph{The Theory of Partitions},
Addison-Wesley, Reading MA, 1976.

\refitem{[C1]}
{I.\,Cherednik},
\emph{Special bases of irreducible representations of a degenerate
affine Hecke algebra},
{Funct.\,Anal.\,Appl.}
{\bf 20}
(1986),
76--78.

\refitem{[C2]}
{I.\,Cherednik},
\emph{A new interpretation of Gelfand-Zetlin bases},
{Duke Math.\,J.}
{\bf 54}
(1987),
563--577.

\refitem{[CP]}
{V.\,Chari and A.\,Pressley},
\emph{Fundamental representations of Yangians and
singularities of R-matrices},
{J.\,Reine Angew.\,Math.}
{\bf 417}
(1991),
87--128.

\refitem{[D1]}
{V. Drinfeld},
\emph{Hopf algebras and the quantum Yang-Baxter equation},
{\text{Soviet} Math.\ Dokl.}
{\bf 32}
(1985),
254--258.

\refitem{[D2]}
{V. Drinfeld},
\emph{A new realization of Yangians and quantized affine algebras},
{\text{Soviet} Math.\ Dokl.}
{\bf 36}
(1988),
212--216.

\refitem{[FM]}
{E.\,Frenkel and E.\,Mukhin},
\emph{Combinatorics of q-characters of finite-dimensional representations
of quantum affine algebras}, preprint
{\tt math.QA/9911112}.

\refitem{[K]}
{M.\,Kashiwara},
\emph{On level zero representations of quantized affine algebras}, preprint 
{\tt math.QA/0010293}.

\refitem{[KMT]}
{N.\,Kitanine, J.-M.\,Maillet and V.\,Terras},
\emph{Form factors of the XXZ Heisenberg spin-$\]\frac12$ finite chain},
{Nuclear Phys.}
{\bf B\,554}
(1999),
647--678.

\refitem{[LNT]}
{B.\,Leclerc, M.\,Nazarov and J.-\ns Y.\,Thibon},
\emph{Induced representations of affine Hecke algebras
and canonical bases of quantum groups}, preprint
{\tt math.QA/0011074}.

\refitem{[M]}
{A.\,Molev},
\emph{Irreducibility criterion for tensor products
of Yangian evaluation modules}, preprint
{\tt math.QA/0009183}.

\refitem{[MNO]}
{A.\,Molev, M.\,Nazarov and G.\,Olshanski},
\emph{Yangians and classical Lie algebras},
{Russian Math.\ Surveys}
{\bf 51}
(1996),
205--282.

\refitem{[MT]}
{J.-M.\,Maillet and V.\,Terras},
\emph{On the quantum inverse scattering problem},
{Nuclear Phys.}
{\bf B\,575}
(2000),
627--644.

\refitem{[N]}
{M.\,Nazarov},
\emph{Yangians and Capelli identities},
{Amer.\,Math.\,Soc.\,Transl.}
{\bf 181}
(1997),
139--164.

\refitem{[NO]}
{M.\,Nazarov and  G.\,Olshanski},
\emph{Bethe subalgebras in twisted Yangians},
{Comm. Math.\,Phys.}
{\bf 178}
(1996),
483--506.

\refitem{[NT1]}
{M.\,Nazarov and V.\,Tarasov},
\emph{Representations of Yangians with Gelfand-Zetlin bases},
{J.\,Reine Angew.\,Math.}
{\bf 496}
(1998),
181--212.

\refitem{[NT2]}
{M.\,Nazarov and V.\,Tarasov},
\emph{On irreducibility of tensor products of Yangian modules},
{Internat.\,Math.\,Research Notices}
(1998),
125--150.

\refitem{[O]}
{G.\,Olshanski},
\emph{Representations of infinite-dimensional classical
groups, limits of enveloping algebras, and Yangians},
{Adv.\,Soviet\,Math.}
{\bf 2}
(1991),
1--66.

\refitem{[T1]}
{V.\,Tarasov},
\emph{Structure of quantum L-operators for the R-matrix of the XXZ-model},
{Theor.\,Math.\,Phys.}
{\bf 61}
(1984),
1065--1071.

\refitem{[T2]}
{V. Tarasov},
\emph{Irreducible monodromy matrices for the R-matrix of the XXZ-model
and local lattice quantum Hamiltonians},
{Theor.\,Math.\,Phys.}
{\bf 63}
(1985),
440--454.

\bye